\documentclass[12pt]{book}

\usepackage[english]{babel}
\usepackage{times,amsfonts,amsmath,amssymb,dsfont,mathrsfs,amsthm} 
\usepackage{pdfsync}
\usepackage{graphicx} 
\usepackage{mathabx}
\usepackage{cases}
\usepackage{caption}
\usepackage{subcaption}
\usepackage{enumitem}

\usepackage{color}

\definecolor{gr}{rgb}   {0.,   0.69,   0.23 }
\definecolor{bl}{rgb}   {0.,   0.5,   1. }
\definecolor{mg}{rgb}   {0.85,  0.,    0.85}
\definecolor{or}{rgb}   {0.9,  0.5,   0.}

\definecolor{webred}{rgb}{0.75,0,0}
\definecolor{webgreen}{rgb}{0,0.75,0}
\usepackage[citecolor=webgreen,colorlinks=true,linkcolor=webred]{hyperref}

\newtheorem{theorem}{Theorem}[chapter]
\newtheorem{proposition}[theorem]{Proposition}
\newtheorem{lemma}[theorem]{Lemma}
\newtheorem{corollary}[theorem]{Corollary}

\theoremstyle{definition}
\newtheorem{definition}[theorem]{Definition}
\newtheorem{notation}[theorem]{Notation}

\theoremstyle{remark}
\newtheorem{remark}[theorem]{Remark}

\newcommand{\beq}{\begin{equation}}
\newcommand{\eeq}{\end{equation}}

\newcommand{\Bk}{\color{black}}

\newcommand{\N}{\mathbb{N}}
\newcommand{\R}{\mathbb{R}}
\newcommand{\C}{\mathbb{C}}
\newcommand{\Z}{\mathbb{Z}}

\newcommand{\rd}{{\mathrm d}}
\newcommand{\one}{\mathds{1}}

\newcommand{\dS}{{\mathbb{S}}}

\newcommand{\cA}{\mathcal{A}}
\newcommand{\cB}{\mathcal{B}}

\newcommand{\cE}{\mathcal{E}}
\newcommand{\cF}{\mathcal{F}}
\newcommand{\cL}{\mathcal{L}}
\newcommand{\cI}{\mathcal{I}}
\newcommand{\cN}{\mathcal{N}}
\newcommand{\cO}{\mathcal{O}}
\newcommand{\cP}{\mathcal{P}}
\newcommand{\cQ}{\mathcal{Q}}

\newcommand{\cS}{\mathcal{S}}
\newcommand{\cT}{\mathcal{T}}
\newcommand{\cU}{\mathcal{U}}
\newcommand{\cV}{\mathcal{V}}

\newcommand{\cX}{\mathcal{X}}

\newcommand{\gh}{\mathfrak{h}}

\newcommand{\gA}{\mathfrak{A}}
\newcommand{\gC}{\mathfrak{C}}
\newcommand{\gD}{\mathfrak{D}}

\newcommand{\gP}{\mathfrak{P}}
\newcommand{\gS}{\mathfrak{S}}

\newcommand{\ogP}{\overline{\mathfrak{P}}}

\newcommand{\bel}{\begin{equation} \label}
\newcommand{\ee}{\end{equation}}

\newcommand{\dD}{\mathbb{D}}
\newcommand{\dP}{\mathbb{P}}

\newcommand{\dx}{\mathbb{X}}

\newcommand{\ogD}{\overline\gD}

\newcommand{\sE}{\mathscr{E}}

\newcommand\eps{\epsilon}
\newcommand\spec{\gS}
\newcommand\eff{\operatorname{eff}}
\newcommand\ran{\operatorname{ran}}

\newcommand\dom{\operatorname{Dom}}
\newcommand\supp{\operatorname{supp}}
\newcommand\curl{\operatorname{curl}}

\newcommand\Id{\operatorname{\mathbb{I}}}
\renewcommand\Re{\operatorname{Re}}
\newcommand\Tr{\operatorname{Tr}}

\makeatletter

\begin{document}

\title{Magnetic fields and boundary conditions in spectral and asymptotic analysis}

\author{Nicolas Popoff}
%

\date{\today}

\maketitle
\tableofcontents

\chapter{Introduction}

This memoir is devoted to a part of the results from the author about two topics: in the first part, the asymptotics of the low-lying eigenvalues of Schr\"odinger operators in domains that may have corners, and in the second part, the analysis of the thresholds of a class of fibered operators. The main common object is the magnetic Laplacian, and the two parts are connected through the study of model problems in unbounded domains. In this short introduction, we present briefly our concerns, without entering the quantitative details.

\section{Low-lying eigenvalues in corner domains}
\label{S:intro1}
In this section we present our problematics around the asymptotics of the low-lying eigenvalues of Schr\"odinger operators in corner domains. More precise definitions and references will be found in Chapters \ref{C:corner}-\ref{C:magnetic}.

In Part \ref{P:lls}, we will mainly consider two operators, acting on functions of a domains $\Omega \subset \R^{n}$ having a boundary: the Laplacian with a Robin type boundary condition $h^{1/2}\partial_{\nu}u= u$, where $h>0$ and $\partial_{\nu}$ denotes the outward normal derivative; and the magnetic Laplacian $(-ih\nabla-A)^2$, completed with natural Neumann boundary conditions, where $A$ is a magnetic potential associated with a given magnetic field $B$. Our analysis extends to other operators, such that the Robin Laplacian with a variable coefficient, or the Schr\"odinger operator with a $\delta$-interaction supported by a hypersurface with corners. In this introduction, we use the generic notation $L_{h}$ for one of these operators. 
In our framework, when considered on a certains class of bounded domains, called {\it corner domains}, these operators are self-adjoint with compact resolvent. We denote by $(\lambda_{j}(h))_{j\geq1}$ their (increasing sequence of) eigenvalues, and our concern is to determine the asymptotics of these eigenvalues in the semi-classical limit $h\to0$, a problematic motivated by physical models such as those involved in the theory of surface superconductivity (\cite{FouHel10}). We will mainly focus on $j=1$, although some of the results will also hold for higher eigenvalues.

For a semi-classical Schr\"odinger operator $-h^2\Delta+V$, with a bounded from below and confining potential $V$, the semi-classical limit of the low-lying eigenvalues is determined at the first order by the the infimum of the potential, and the associated eigenfunctions are localized near the minimum (\cite{HeSj84,DiSj99}). In our case, the semi-classical limit of the eigenvalues will be driven by the geometry (of the domain, and of the magnetic field for the magnetic Laplacian). But it is not clear how to define a quantity which will play the same role as the minimum of $V$ in the standard case does. 

The operators considered have an homogeneity property with respect to dilations. Indeed, if the operator $L_{h}$ is considered on a cone, with its coefficients frozen in a suitable way, then it is unitarily equivalent to $hL_{1}$. 
This property leads to a natural idea: for domains with corners (see Section \ref{S:presentationcorner} for a rigorous definition), one will have to look at the operator on the tangent geometries. Assume that $\Omega$ is domain with corners, given a point $x\in \overline{\Omega}$, firstly, we need to define a tangent operator, with coefficients frozen in a suitable way, with $h=1$, on the tangent cone at $x$. Denote by $E(x)$ the bottom of the spectrum of this operator. Then, the general expected result, which has appeared to be true for all the particular cases treated in the literature, is that 
\bel{E:megameta}
\lambda_{1}(h)\underset{h\to0}{\sim} h\sE, 
\ee
with 
\bel{D:lle}
\sE=\sE(\Omega)=\inf_{x\in\overline{\Omega}}E(x).
\ee
We will give more references later, we refer to those in the book \cite{FouHel10} and \cite{BoDauPof14} for the magnetic Laplacian, and to \cite{LevPar08} for the Robin Laplacian. The first order term in the asymptotics is therefore linear with respect to $h$. This comes from the homogeneity property described above. The minimization can be understood by keeping in mind that the eigenfunctions associated with $\lambda_{1}(h)$ will tend to concentrate near some point $\overline{\Omega}$, in particular near points where $E(\cdot)$ is minimum. We have called the function $E$ the {\it local energy}.

To prove such an asymptotics for general domains rises several problem: 
\begin{itemize}
\item In the case treated, the local energy is often discontinuous when changing of strata inside of a corner domain. Therefore, it is important to show that the minimum of the local energy is reached, and is non-degenerate, in the sense that $\sE\notin \{0,-\infty\}$. Moreover, is it possible to determine the geometry minimizing $E$?
\item Is it possible to show \eqref{E:megameta}, together with an estimate for $\lambda_{1}(h)-h\sE$, without additional hypotheses?
\item Is it possible to have a more precise asymptotics of $\lambda_{1}(h)-h\sE$, together with an asymptotics of the higher eigenvalues, under additional hypotheses on the geometry?
\end{itemize}
The answer to the first question, mainly developed in \cite{BoDauPof14}, is to use the {\it singular chains} of a corner domain, recursively defined, from \cite{Dau88}. The set of singular chains of $\Omega$ extends the points of $\overline{\Omega}$, in the sense that it takes into account the local geometry of $\Omega$. The idea behind this procedure is to desingularize the domain, it originates from \cite{Kon67}, and is not far from the concept of iterated blow-up (\cite{Mel15,Gri17}).
We will consider the local energy on singular chains, and show that it is lower semi-continuous, and therefore reaches its infimum. Moreover, we will show a monotonicity property which, roughly, expresses as ``if two geometries can be compared, then the more singular one will have the lowest local energy''. The model example is the one of the magnetic Laplacian in a wedge and of its two faces, studied in \cite{Pof15}.

We have succeeded in the second question, by giving separately a lower bound and an upper bound for $\lambda_{1}(h)-h\sE$. The lower bound relies on a classical idea: a suitable partition of the unity, together with the well known IMS formula, should allow to compare the operator to local models. However this procedure cannot be done directly, due to the possible blow-up of one of the principal curvatures in corner domains, indeed large principal curvatures will result in a bad estimate when approximating the operator by a frozen operator on its tangent geometry. This problem is solved by a multiscale analysis, adapted to the recursive definition of a corner domain. The upper bound relies on the construction of a test-function, whose energy is close to $h\sE$. This test-function will come from a tangent geometry in which the model operator has an eigenfunction. The existence of such a tangent geometry is not an easy question, and is linked to the first question. In particular, the local energy on the tangent {\it substructures} of a cone (such as the faces for a wedge) will play the role of a threshold in the spectrum of the tangent operator. Notice that the construction of the test functions involves also a multiscale procedure in order to counterbalance the possible blow up of the principal curvatures. 

The third question will need more hypotheses on the structure of the local energy near its minimum. If one thinks of the local energy as an effective potential, this is coherent with the harmonic approximation, in which the non-degeneracy of the minimum of the potential provides an asymptotic expansion in powers of the semiclassical parameter (\cite{DiSj99}). But in our case, it is not obvious what are the natural hypotheses on the geometry, since the local energy does not have an explicit expression. More or less, one may think of three kinds of natural {\it generic }setting:
\begin{itemize}
\item Corner concentration. This is the case when the local energy has a discontinuous jump at its minimum. The model case is the one of a polygonal domain, in which the local energy is minimum at a corner, and has a gap with its values on the sides. In some sense, the local problem at this corner will dominate all the others, concerning the asymptotics of the first eigenvalue. This case has been treated under some assumptions in \cite{Bon06,BonDau06} for the magnetic Laplacian and \cite{HePank15,Kha18} for the Robin Laplacian. 
\item Wells of the local energy. This is the case when the local energy reaches its infimum continuously inside a stratum. The model case is the one of a varying magnetic field, whose intensity is minimum at some given point in the interior of a stratum, for example at some given point of the boundary of a regular domain. We cite mainly \cite{LuPan99-2,HeMo01,FouHel10,Ray13,Ray3d09}. A domain with an edge whose opening has a non-degenerate extrema may enter also this framework (\cite{Popoff,PofRay13}). Most of the time, a generic assumptions on the geometry has to be added, leading to the non-degeneracy of the local energy. 
\item Submanifold wells. This is the case when the local energy is minimal on a submanifold. Typical cases are the Robin Laplacian in a bounded regular domain (the local energy is constant on the whole boundary) (\cite{Pank13,ExMinPar13,PankPof15,HeKa15}), and the magnetic Laplacian in a regular domain with a constant magnetic field (\cite{HeMo01,HeMo04,FouHe06}). 
\end{itemize}

The first two cases are essentially local, in the sense that the behavior of the local energy at one point will determine the next terms in the asymptotics, at least for those which are powers of $h$. Note however that for a straight polyhedron, the remainder may be exponentiall small, and the behavior of $\lambda_{1}(h)-h\sE$ depend on the global geometry, and is given by  tunneling effect in the case of symmetry.

The last case poses the question of the existence of an effective Hamiltonian, defined on the submanifold wells, leading the asymptotics of the low-lying eigenvalues. This problem is solved for the Robin Laplacian in \cite{PankPof16}, where we have introduced a semi-classical Laplace operator on the boundary, involving the mean curvature. Still, the existence of such an effective operator is still a challenging question for the Neumann magnetic Laplacian, even in dimension two, and would be an important step toward the understanding of the tunneling effect for the magnetic Laplacian in symmetrical regular domains (\cite{BoHerRay16}).

In Chapter \ref{C:corner}, we define the corner domains of a Riemannian manifold from \cite{Dau88}, together with their singular chains. We present the analysis from \cite{BoDauPof14}. We also present a formal IMS formula, providing a mechanism for a lower bound of the first eigenvalue of our operators in such domains, including error terms which will be detailed later, depending on the context.
In Chapter \ref{C:robin}, we introduce the Robin Laplacian in such a domain. We present the analysis of the operator on the tangent cones, and the recursive procedure leading to a two-side estimate for the first eigenvalue in \cite{BruPof16}. For regular domains, a more precised asymptotics was obtained in \cite{PankPof15,PankPof16}, showing the existence of an effective Hamiltonian, defined on the boundary. We apply our asymptotics to a reverse Faber-Krahn inequality. We also present the analysis for a Robin Laplacian with a vanishing coefficient, a case where the operator is non-self adjoint (\cite{NazPof18}).

In Chapter \ref{C:magnetic}, we introduce the Laplacian with magnetic field and Neumann boundary condition. This operator is, in some sense, more intricate to analyse, because the magnetic field combines with the geometry in the determination of the minimum of the local energy. We present the results analogous to those of chapter \ref{C:robin} for corner domains, coming from \cite{Pof15,BoDauPof14}, and enlighten the differences in the treatment of the problems. Mainly, a recursive analysis does not seem available and we proceed to an exhaustion of model problems to reach the asymptotics of the first eigenvalue in dimension 3.  We also give improvements of the asymptotics under stronger hypotheses, from \cite{PofRay13}.

Notice that the results from Chapter \ref{C:magnetic} are anterior to those of Section \ref{SS:Robincoin}. 

\section{The thresholds of translationally invariant magnetic Laplacians}

In a second part, we consider Schr\"odinger operators with translationally invariant magnetic fields. Our study of these Hamiltonians is mainly motivated by two different contexts: firstly, they appear as local model problems in the study of the semi-classical magnetic Laplacian, in particular the analysis of the Laplacian with a constant magnetic field in a half-plane is a necessary step in the semi-classical asymptotics of the Laplacian in a regular domain. Secondly, the associated quantum systems present interesting transport properties in the direction of invariance, and these models are used in the understanding of the Quantum Hall Effect.

Our operators are of the form
\bel{E:translationnally} 
H_{0}=(-i\nabla -A)^2, 
\ee
acting on $L^2(\Omega)$, where $\Omega:=\omega\times \R$, and $\omega\subset \R^{n-1}$ (with $n=2$ or $n=3$), moreover, the magnetic feld $B=\curl A$ does not depend on the last variable. Our three main models are Schr\"odinger operators in a half-plane with a constant magnetic field (\cite{dBiePu99,HeMo01,BruMirRai13}), in $\R^{2}$ with a magnetic field varying in only one direction (sometimes called the {\it Iwatsuka model}, \cite{Iwa85,ManPur97,HisSoc14}), and in $\R^{3}$ with a cylindrical-invariant magnetic field (\cite{Rai08,Yaf08}), but our approach can be extended to a wide class of other Hamiltonians. 
 
If we denote by $\cF$ the partial Fourier transform along the direction of invariance, our operator enters the framework of operators which can be fibered over a real analytic manifold $M$:
\bel{E:fibering}
\cF H_{0} \cF^{*}= \int_{k\in M}^{\bigoplus}\gh (k) \rd k,
\ee
with $M=\R$ in our case, and where $(\gh(k))_{k\in M}$ is a family of positive self-adjoint operators operators (they are of Sturm-Liouville type when $\omega\subset \R$). Most of the time, $\gh(k)$ has compact resolvent and we denote by $\lambda_{j}(k)$ its eigenvalues. Their are called the {\it band functions}, or {\it dispersion curves} of the system. The spectrum of $H_{0}$ is 
\bel{E:fiberspectrum}
 \sigma(H_{0})=\bigcup_{j\geq1}\overline{\lambda_{j}(\R)},
 \ee
it is absolutely continuous, provided that the band functions are not constant. Moreover, some energies in the spectrum enjoy remarkable properties, they are thresholds in the spectrum. A general theory for thresholds of analytically fibered operators exists, \cite{GeNi98}, but one of the properties of the class described in that article is that the band functions are proper. This is not the case in our magnetic models, since the band functions may tend to finite limits as $k\to\infty$, giving rise to a new kind of thresholds. 

Our goal is not to develop a general theory for these thresholds, but to illustrate several phenomena typical of Hamiltonians whose band functions tend to finite limit. Unlike to critical points of band function, there ``is no Taylor expansion'' near these values, and a first step is to provide an asymptotics expansion of $\lambda_{j}$ at infinity. This is done by standard tools coming from the harmonic approximation, the operators $\gh(k)$ turning to be of semi-classical type as $k\to\infty$ (\cite{Popoff,HisPerPofRay16}). It is known that quantum states localized in energy far from the thresholds enjoy transport and localization properties in the following sense: they bear a current which can be bounded from below, and they are usually called {\it edge states}, because they are small far from the boundary in Quantum Hall systems (\cite{dBiePu99,DomGerRai11}). On the opposite side, as we explain in Section \ref{S:current}, the quantum states localized near these thresholds have a component with very small velocity along the invariance direction, moreover they are localized at infinity, both in space and in frequency. This analysis can be found in \cite{HisPofSoc14,MirPof18}.

Then, we are interested in suitable perturbations of these operators. The essential spectrum will remain the same, but eigenvalues may appear. To evaluate how many discrete eigenvalues are created is a widely studied question. In particular, these eigenvalues tend to accumulate near the edges of the essential spectrum, that are the end-points of the set \eqref{E:fiberspectrum}. In the case of a constant magnetic field in the whole space, the question is well known, see \cite{AvrHerSi78}, \cite{RaiWar02} (and the references therein). When such an edge corresponds to a critical point of a band function, through a localization in frequencies, it appears that an effective Hamiltonian governs the asymptotics of the eigenvalues in this zone, see \cite{Rai92} for the analysis in Schr\"odinger operators with periodic potentials. But near a threshold which is a limit of a band function, this procedure does not work directly. An effective operator is given in \cite{BruMirRai13} for a constant magnetic field in a half-plane with Dirichlet conditions.
The analogous of these questions inside the essential spectrum is to give the behavior of the spectral shift function associated with the perturbation. It is expected that this function may be singular at the thresholds. In Section \ref{S:perturbIwa}, we present the results from \cite{MirPof18}: we consider the Iwatsuka model submitted to an electric perturbation, and we provide the {\it a priori} and the precise behavior of this function near thresholds, depending on the hypotheses on the decay of the magnetic field and the electric perturbation. 

In Section \ref{S:magneticwire}, we present the Schr\"odinger operator with a magnetic field created by a infinite rectilinear wire, already considered in \cite{Yaf03}. This model possesses an additional difficulty: The bottom of the spectrum corresponds to an accumulation of band functions, see Figures \ref{F1}--\ref{F2}. We study conditions on the electric perturbation for having the finiteness of eigenvalues below the essential spectrum (\cite{BruPof15}).

Numerous questions remains unsolved for such operators, both in particular cases or in a more global approach. For the model case of the Dirichlet half-plane with constant magnetic field, perturbed by a compactly supported potential, the precise behavior of the number of eigenvalues is still not known.
Moreover, for these models, the nature of these thresholds, as branching points of the resolvent, seems to be an interesting question. We expect that these branch points will have an original behavior, compared to the branch points of the Laplacian with a periodic potential, well described (\cite{Ge90}). This would be a starting point in order to tackle the analysis of the resonances near a threshold for perturbations of $H_{0}$. 

\part{Low lying spectrum of Laplacians in corner domains}
\label{P:lls}
\chapter{Operators in corner domains}
\label{C:corner}
In this chapter we present the class of domains $\Omega$ with corners and their singular chains $\dx$, extending the points of $\overline{\Omega}$, in the spirit of \cite{Dau88}. We introduce the local energy of an operator $L_{h}$, $E(\dx)$ as the infimum of the spectrum of the tangent operator associated with the chain $\dx$. In the last section, we present a general IMS formula, based on a multiscale analysis, providing a lower bound for the first eigenvalue of $L_{h}$, as $h\to0$. At this stage, the operator $L_{h}$ is not specified, and the method will be applied to particular cases in Chapter \ref{C:robin} and \ref{C:magnetic}.


\section{Presentation of corner domains}
\label{S:presentationcorner}
The operators we will consider share the property to have a nice homogeneity property with respect to dilations. If one thinks that semi-classical asymptotics requires the analysis of {\it local models}, being defined as the frozen operator on the tangent geometry at a point, it is expected to consider domains which locally are close to a cone, which is defined as a subset of $\R^{n}$ invariant by positive dilations. In this spirit, the main class of domains in which we will work are the {\it corner domains}, defined in \cite{Dau88}, in the spirit of \cite{MazPla77}:

\begin{definition}[\sc Class of corner domains]\index{Corner domain|textbf}
\label{D:cornerdomains}
The classes of corner domains $\gD(M)$\index{Corner domain!class of --\quad$\gD(M)$\ |textbf} ($M=\R^n$ or $M=\dS^n$) and tangent open cones $\gP_{n}$ are defined as follows: 

\noindent{\sc Initialization}, $n=0$: 
\begin{enumerate}
\item $\gP_0$ has one element, $\{0\}$,
\item $\gD(\dS^0)$ is formed by all (non empty) subsets of $\dS^0$.
\end{enumerate}

\noindent{\sc Recurrence}: For $n\ge1$,
\begin{enumerate}
\item $\Pi\in\gP_n$ if and only if its section of $\Pi\cap \dS^{n-1}$ belongs to $\gD(\dS^{n-1})$,
\item \label{machin}
$\Omega\in\gD(M)$ if and only if for any $x_{0}\in\overline\Omega$, 
there exists an open tagent cone $\Pi_{x_{0}}\in\gP_n$ to $\Omega$ at $x_{0}$. 
\end{enumerate}
\end{definition}
The existence of a tangent cone is linked to a diffeomorphism $\psi_{x_{0}}: \cU \to \cV$, where $\cU$ (respectively $\cV$) is a neighborhood of $x_{0}$ (respectively of 0), and such that $\psi_{x_{0}}(\cU\cap \Omega)=\cV \cap \Pi_{x_{0}}$ and $\psi_{x_{0}}(\cU\cap \partial\Omega)=\cV \cap \partial\Pi_{x_{0}}$. The open set $\cU$ is called a map-neighborhood of $x_{0}$.

Note that $\gD(\R^{n})$ includes smooth domains. Let us introduce a subclass of corner domains. 
\begin{definition}[\sc Class of polyhedral cones and domains] \index{Polyhedral domain|textbf} \index{Polyhedral cone}
The classes of polyhedral domains $\ogD(M)$\index{Polyhedral domain!set of -- \quad$\ogD(M)$\ |textbf} ($M=\R^n$ or $M=\dS^n$) and polyhedral cones $\ogP_{n}$\index{Polyhedral cone!set of -- \quad$\ogP_{n}$\ |textbf} are defined as follows: 

\begin{enumerate}
\item The cone $\Pi\in \gP_n$ is a polyhedral cone if its boundary is contained in a finite union of subspaces of codimension $1$. We write $\Pi\in\ogP_n$.
\item  The domain $\Omega\in\gD(M)$ is a polyhedral domain if all its tangent cones $\Pi_x$ are polyhedral. We write $\Omega\in\ogD(M)$.
\end{enumerate}
\end{definition}
Roughly, one may think that the regular part of a polyhedral cone has zero curvature. On the contrary, every cone in $\gP_n\setminus \ogP_{n}$ has an unbounded principal curvature near the origin, by a direct dilation argument. As a consequence, the polyhedral corner domains have bounded curvatures, but the non polyhedral have not. 

\begin{figure}[ht]
\noindent\centering
    \begin{subfigure}[h]{0.32\textwidth}
\includegraphics[keepaspectratio=true,width=5.3cm]{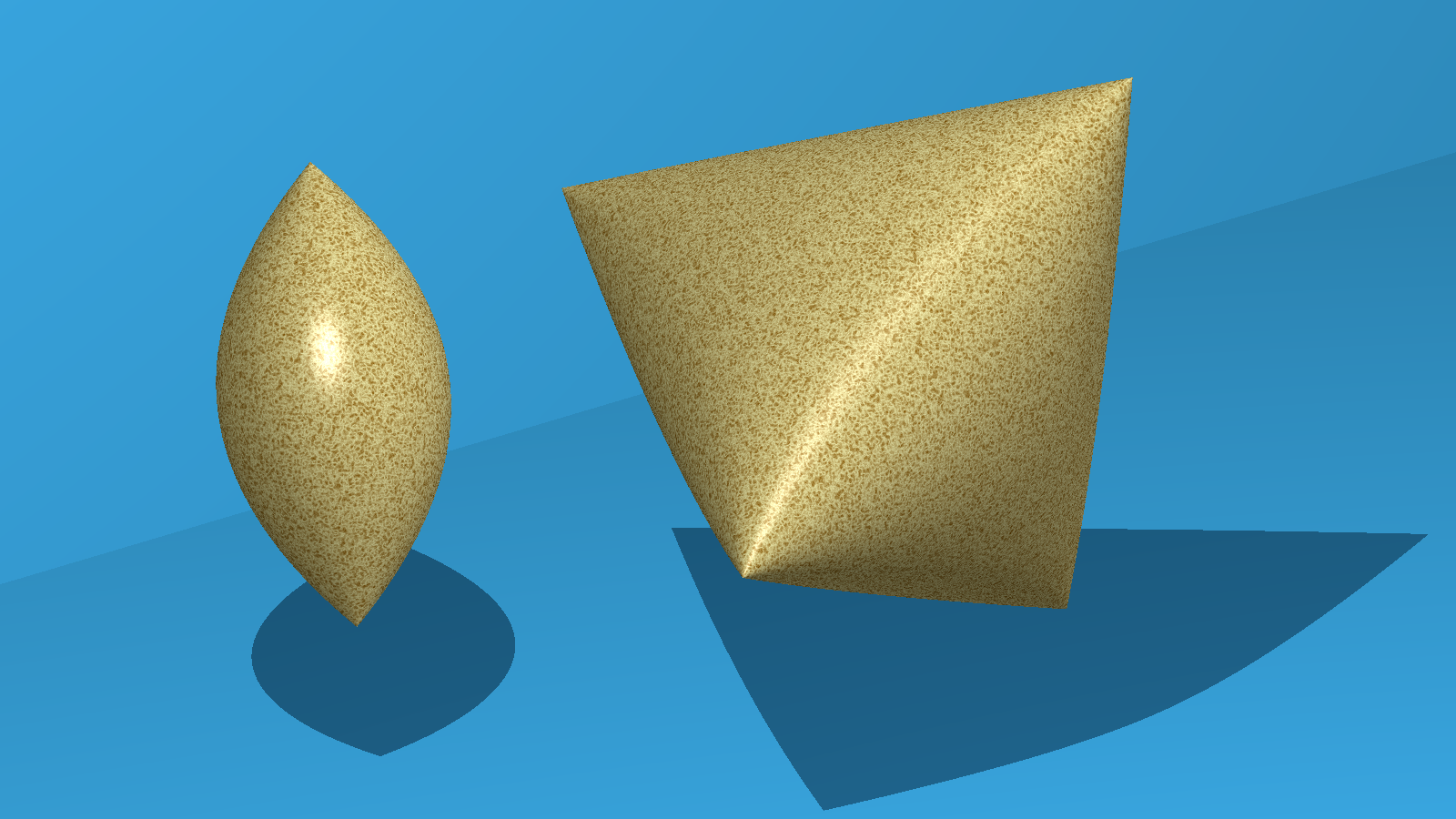}
        \caption{Domains with corners}
        \label{F1A}
    \end{subfigure}
\
            \begin{subfigure}[h]{0.32\textwidth}
\includegraphics[keepaspectratio=true,width=5.3cm]{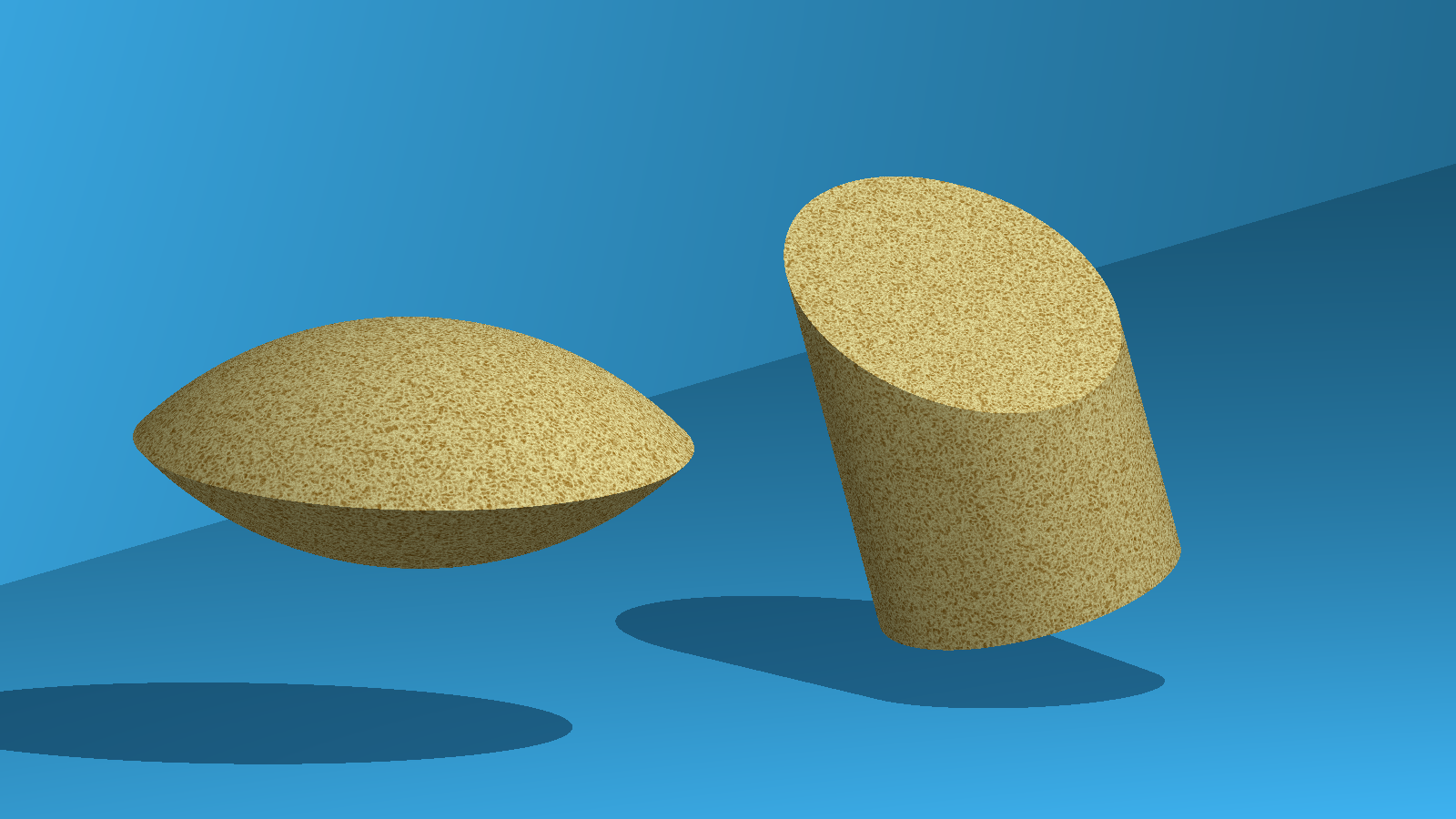}
        \caption{Domains with edges}
        \label{F1B}
    \end{subfigure}
    \
        \begin{subfigure}[h]{0.32\textwidth}
\includegraphics[keepaspectratio=true,width=5.3cm]{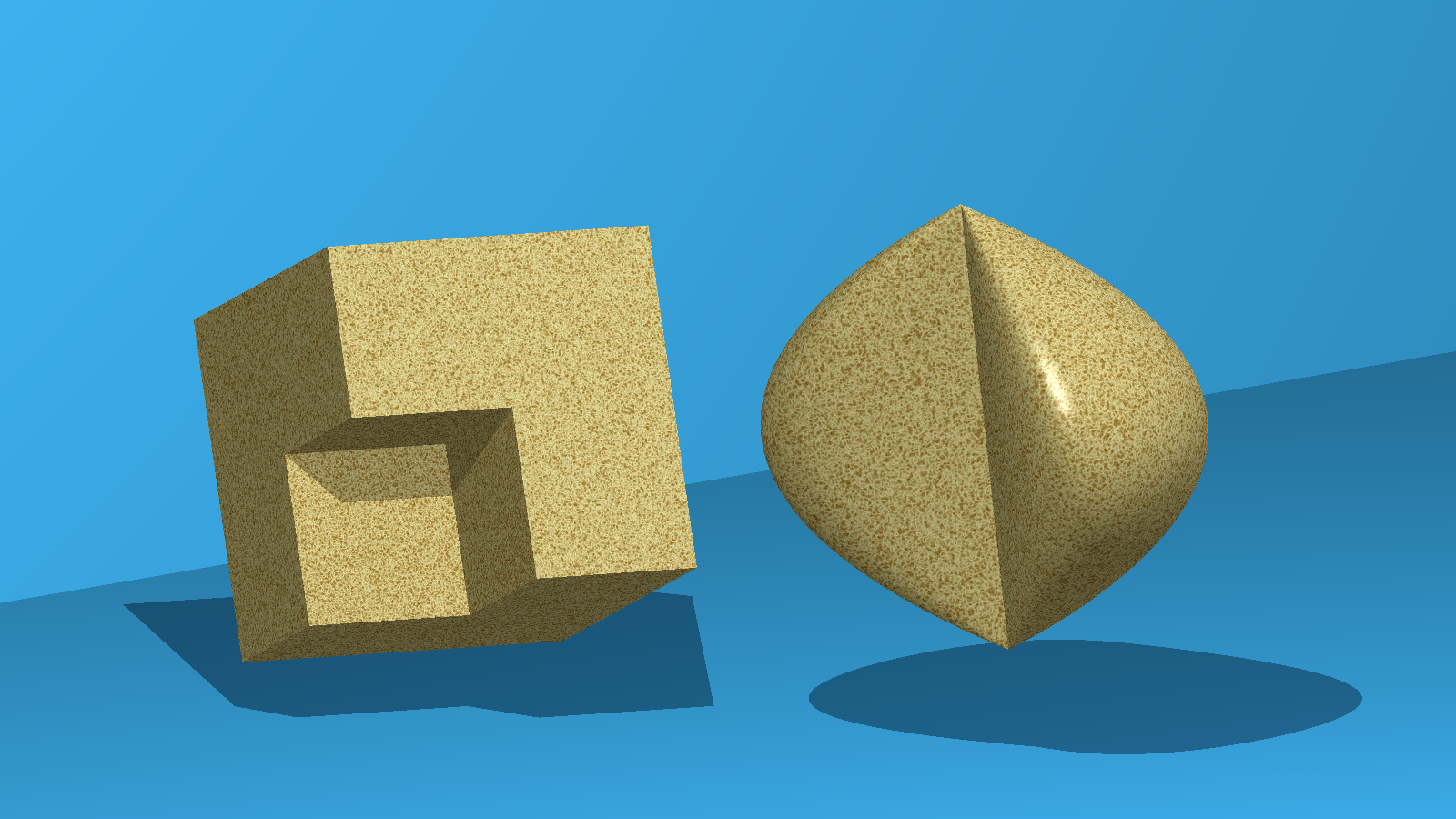}
        \caption{Domains with both}
        \label{F1C}
    \end{subfigure}
\caption{Examples of 3D corner domains (Figures drawn by M. Costabel with POV-ray)}
\label{Fcorner}

\end{figure}

In dimension 2 the elements of $\gP_2$ are $\R^2$ and all plane sectors with opening $\alpha\in(0,2\pi)$, denoted by $S_{\alpha}$, including half-planes ($\alpha=\pi$).  
Therefore, the elements of $\gD(\R^2)$ are the regular domains, and the curvilinear polygons with piecewise non-tangent smooth sides (corner angles $\alpha\neq0,\pi,2\pi$). In particular, in dimension 2, we have $\gP_2=\ogP_2$ and $\gD(\R^2)=\ogD(\R^{2})$. This is not true in dimension 3, as shows the example of a circular cone, which is not in $\ogP_{3}$. In the examples of figure \ref{Fcorner}, the corner domains of Figure \ref{F1A} are not polyhedral, whereas those of Figure \ref{F1B} are. 

 A cone $\Pi\in \gP_{n}$ being given (up to a rotation) in the form $\R^{\nu}\times \Gamma$, with $\nu\geq0$ maximal for such a form, we denote by $d(\Pi):=n-\nu$ the {\it irreducible dimension} of $\Pi$. Now, for a corner domain $\Omega \in \ogD(M)$, we denote by 
\begin{equation}
\label{eq:Ad}
   \gA_d(\Omega) = \{x\in\overline\Omega,\quad d(\Pi_{x})=d\}.
\end{equation}
In \cite{BoDauPof14}, we prove that the corner domains admit a stratification:
\begin{proposition}
Let $\Omega\in \gD(M)$ and let $0\leq d \leq n$. Then the connected component of $\gA_d(\Omega)$ are submanifold of codimension $d$.
\end{proposition}
This proposition illustrates the local structure of corner domains. In particular, $\gA_{0}(\Omega)$ is the interior of $\Omega$, $\gA_{1}(\Omega)$ is the regular part of the boundary, $\gA_{2}(\Omega)$ are the edges (when $n\geq3$), and $\gA_{n}(\Omega)$ are the vertices, which are therefore isolated. 


\section{Operators on singular chains and local energy}
\label{S:cornerpure}

In section we present the geometrical setting, together with Theorem \ref{T:Fsci}, that will help to determine the first order term in the lowest eigenvalues of the operator considered.



An integer $p\geq0$ being given, a singular chain $\dx=(x_{0},\ldots,x_{p})$ is a sequence of points defined recursively as follow: $x_{0}$ is set in $\overline{\Omega}$, then denote by $\Pi_{x_{0}}=U_{0}(\R^{\nu_{0}}\times \Gamma_{x_{0}})$ the decomposition of the tangent cone at $x_{0}$, where $U_{0}$ is a rotation. Then $x_{1}$ is picked in $\overline{\Omega_{x_{0}}}$, where $\Omega_{x_{0}}:=\Gamma_{x_{0}}\cap \dS^{d_{0}-1}$. Notice that $\Omega_{x_{0}}\in \gD(\dS^{d_{0}-1})$, therefore there exists a tangent cone at $x_{1}$. The other points $(x_{i})_{i\geq2}$ are defined in the same way, recursively.

To this singular chain is associated a cone $\Pi_{\dx}$ as follows: $\Pi_{(x_{0})}=\Pi_{x_{0}}$, $\Pi_{(x_{0},x_{1})}=U_{0}(\R^{n-d_{0}}\times \langle x_{1}\rangle \times \Pi_{x_{1}})$, where $\langle x_{1}\rangle$ is the vector space generated by $x_{1}$, (recall that $\Pi_{x_{1}}\subset \R^{d_{0}-1}$ denotes the tangent cone at $x_{1}\in \overline{\Omega_{x_{0}}}$), and so on for a chain $\dx=(x_{0},\ldots,x_{p})$. These cones are called the tangent substructures of $\Omega$ at $x_{0}$. We denote by $\gC(\Omega)$ the set of all singular chains of $\Omega$, and $\gC_{x_{0}}(\Omega)$ the set of chains originating at a given point $x_{0}\in \overline{\Omega}$.

As an illustration, we give below an exhaustion of the chains originating at a vertex $x_{0}\in \gA_3(\Omega)$ of $\overline{\Omega} \in \gD(\R^{3})$: 

There are four possible lengths for chains in $\gC_{x_0}(\Omega)$:
\begin{enumerate}
\item $\dx = (x_0)$ with $\Pi_\dx=\Pi_{x_0}$, the tangent cone (which is not a wedge). It coincides with its reduced cone since $x_0\in \gA_3$. Its section $\Omega_{x_0}$ is a corner domain in  $\dS^2$.
\item $\dx = (x_0,x_1)$ where $x_1\in \overline \Omega_{x_0}$. 
\begin{enumerate}
\item If $x_1$ is interior to $\Omega_{x_0}$, $\Pi_\dx=\Pi_{(x_{0},x_{1})}\equiv \langle x_{1}\rangle\times \R^{2}\equiv \R^3$. No further chain.
\item If $x_1$ is in a side of $\Omega_{x_0}$, $\Pi_\dx\equiv \langle x_{1}\rangle\times \R\times \R^{+}\equiv \R^2\times\R_{+}$ is a half-space.
\item If $x_1$ is a corner of $\Omega_{x_0}$ with opening angle $\theta\in ]0,2\pi[\setminus\{\pi\}$, then we denote by $S_{\theta}$ an open infinite secteur in the plane, and we have, $\Pi_\dx\equiv \langle x_{1}\rangle\times S_{\theta}\equiv \R \times S_{\theta}$ is a wedge. Its edge contains one of the edges of $\Pi_{x_0}$.
\end{enumerate}
\item $\dx=(x_0,x_1,x_2)$ where $x_1\in\partial\Omega_{x_0}$ 
\begin{enumerate}
\item If $x_1$ is in a side of $\Omega_{x_0}$, the reduced tangent cone at $x_{1}$, denoted by $\Gamma_{x_{1}}$, is a half-line: $\Gamma_{x_{1}}\equiv \R^{+}$. In that case, $\Omega_{x_{1}}=\{1\}$, and  $x_2=1$, therefore $\Pi_\dx=\R^3$. No further chain.
\item If $x_1$ is a corner of $\Omega_{x_0}$,  $\Gamma_{x_{1}}$ is a plane sector, its section $\Omega_{x_{1}}$ is an open interval of the unit cercle, and $x_2\in\overline \Omega_{x_{1}}$.
\begin{enumerate}
\item If $x_2$ is an interior point of $\Omega_{x_{1}}$, then $\Pi_\dx=\R^3$. 
\item If $x_2$ is a boundary point of $\Omega_{x_{1}}$, then $\Pi_\dx$ is a half-space.
\end{enumerate}
\end{enumerate}
\item $\dx=(x_0,x_1,x_2,x_3)$ where $x_1$ is a corner of $\Omega_{x_0}$, $x_2\in\partial\cI_{x_0,x_1}$ and $x_3=1$. Then $\Pi_\dx=\R^3$.
\end{enumerate}
Notice that different chains can lead to the same tangent structure. In that case, the chains are called equivalent.  These chains are set with a natural partial order together with a distance: 
\begin{definition}(Order and distance on chains)
\label{D:Ordrechain}
Let $\dx=(x_{0},\ldots,x_{p})$ and $\dx'=(x_{0}',\ldots,x_{p'}')$ be two singular chains in $\gC(\Omega)$. 

We say that $\dx \leq \dx'$ if 
$p\leq p'$ and $x_{j}=x_{j}'$ for all $0 \leq j\leq p$.

 We define the distance\index{Distance of chains|textbf} $\dD(\dx,\dx')\in\R_+\cup\{+\infty\}$ as \index{Distance of chains!$\dD(\dx,\dx')$}
\[
   \dD(\dx,\dx') = \|x_0-x'_0\| + \frac12\bigg\{
   \min_{\substack {L\in\mathsf{BGL}(n) \\ L\Pi_\dx = \Pi_{\dx'}}}
   \|L-\Id_n\|
      +\min_{\substack {L\in\mathsf{BGL}(n) \\ L\Pi_{\dx'} = \Pi_{\dx}}}
   \|L-\Id_n\| \bigg\}\,,
\]
where the second term is set to $+\infty$ if $\Pi_\dx$ and $\Pi_{\dx'}$ do not belong to the same orbit for the action of $\mathsf{BGL}(n)$ on $\gP_{n}$, where $\mathsf{BGL}(n)$ is the semi-group of linear isomorphisms $L$ with norm $\|L\|\le1$.
\end{definition}
In particular, two chains are equivalent if and only if their distance is zero.

Let $\Omega\in \gD(M)$ and $\dx=(x_{0},\ldots,x_{p})\in \gC(\Omega)$, we consider self-adjoint operators bounded from below $L_{h}$, acting on $L^{2}(\Omega)$. In almost all our applications, the form domain of $L_{h}$ will be $H^{1}(\Omega)$. For a given chain $\dx\in \gC(\Omega)$, we need to define the tangent operator on the associated tangent cone, denoted by $L_{h}[\Pi_{\dx}]$, which will be defined precisely according to the context, and is, roughly, the differential operator, frozen on the tangent geometry. We refer to Definitions \ref{N:robin} for the Robin Laplacian, and \ref{N:magnetic} for the magnetic Laplacian. 
%
%

In the cases we will analyze, we will have 
\beq
\label{E:homoge}
L_{h}[\Pi_{\dx}] \equiv h L_{1}[\Pi_{\dx}]=:h L[\Pi_{\dx}]
\ee
We define the local energy on singular chains as
\bel{E:DE}
E(\dx), \ \ \mbox{the bottom of the spectrum of} \ \ L[\Pi_{\dx}].
\ee
When $\dx$ is a chain of length 1, i.e. $\dx=(x)$ with $x\in \overline{\Omega}$, we will simplify by $E(\dx)=E(x)$. We also denote by $\sE(\Omega):=\inf_{\dx \in \gC(\Omega)}E(\dx)$. We will show that $E$ acts as an effective potential, in the sense that its infimum provides the first order in the asymptotics of the first eigenvalue. In particular, the fact that the infimum is reached is an important step stone. We will also show that $\sE$ is involved when determining the essential spectrum of tangent operators.  To proceed in the analysis, we have shown in \cite{BoDauPof14}:

\begin{theorem}
\label{T:Fsci}
Let $F:\gC(\Omega)\to\R$ be continuous and monotonous, for the order and the distance defined in Definition \ref{D:Ordrechain}. Then $F$ is lower semi-continuous. In particular, $F$ restricted to $\overline{\Omega}$ reaches its infimum.
\end{theorem}
We will have to show that $E(\cdot)$, defined as above on singular chains, satisfies the hypotheses of the above theorem.

%

\section{Lower bound in corner domains: IMS formula}
\label{S:Lowerbound}
In this section we present {\it a priori} lower bounds for the first eigenvalue, based on localization formulas when using a partition of the unity. This mechanism works for several Hamiltonian, generically denoted by $L_{h}$ here, which will be presented in the next section.

Let $\cP$ be a set of points of $\overline{\Omega}$, and 
assume that $(\chi_{p})_{p\in \cP}$ is a locally finite, regular, quadratic partition of the unity (i.e. $\sum_{p\in \cP} \chi_{p}^2=1$), depending on 
$h$. Then in the sense of form: 
$$L_{h}=\sum_{p\in \cP}\chi_{p}L_{h}\chi_{p}-h^2\sum_{p\in \cP}|\nabla \chi_{p}|^2.$$
Now, we assume that the the supports of $\chi_{p}$ form a suitable covering of $\overline{\Omega}$, in the sense that the support of each $\chi_{p}$ is included in a map-neighborhood of $p$ and supported in a ball of size $\eps(p)$, with $\eps(p)\to0$ as $h\to0$. Then, still in the sense of form:
$$L_{h}\equiv\sum_{p\in \cP}\chi_{p}L_{h}[\Pi_{p},G_{p}]\chi_{p}-h^2\sum_{p\in \cP}|\nabla \chi_{p}|^2$$
where $L_{h}[\Pi_{p},G_{p}]$ denotes the operator on the tangent cone $\Pi_{p}$ with metric $G_{p}:=J_{p}^{-1}(J_{p}^{-1})^{T}$, $J_{p}$ being the Jacobian of the local diffeomorphism $\psi_{p}$, see Definition \ref{D:cornerdomains}. Now, in regular domains and in polygons, it is standard to approach $G_{p}$ by the identity metric, and to take the homogenous part of the operator.

We denote by $ \kappa(p)$ the $L^\infty$ norm of the curvatures at a point $p$. Error terms, denoted by $R$, depending on the size of $\supp\chi_{p}$, and on $\kappa(p)$, appear. This error term is at least of size $\eps(p)\kappa(p)$. Roughly, using \eqref{E:homoge}, we get (remember that $h \ll 1$):
  \begin{eqnarray}
L=& \sum_{p\in \cP}\chi_{p}hL[\Pi_{p}]\chi_{p}(1+R(\eps(p),\kappa(p))+\cO(h^2\sum_{p\in \cP}\eps(p)^{-2}) \nonumber
\label{E:IMS1}
\\
\geq & \sum_{p\in \cP}\chi_{p}hE(\Pi_{p})\chi_{p}(1+R(\eps(p),\kappa(p)))+\cO(h^2\sum_{p\in \cP}\eps(p)^{-2}) \nonumber
\\
\label{IMSrem} \geq &h\sE(\Omega)+ h\sum_{p\in \cP}R(\eps(p),\kappa(p))+\cO(h^2\sum_{p\in \cP}\eps(p)^{-2})
\end{eqnarray}
 In regular domains, $\kappa$ is bounded, therefore it is enough to take balls of size $\eps(p)=h^{\delta}$, with $\delta<\frac{1}{2}$ chosen in order to optimize the remainders (see \cite{HeMo01} for the magnetic Laplacian).

  Unfortunately, without a refinement, this procedure does not provide a lower bound in general corner domains, because $G_{p}-I$ is linked to the curvatures of the boundary, which may be unbounded. The idea, developed in \cite{BoDauPof14} and extended in \cite{BruPof16}, is to take $\eps(p)$ sufficiently small to counterbalance the possible blow up of the curvatures at $p$. In dimension $n=3$, at a point $p$ near a conical point $x_{0}\in \gA_{0}(\Omega)$ (that is a point such that $\Pi_{x_{0}}\notin \ogP_{3}$), the curvature $\kappa(p)$ can only be controlled by $|p-x_{0}|^{-1}$. In this case, the procedure is as follows: 

\begin{itemize}
\item Take a ball of size $\eps(x_{0})=h^{\delta_{0}}$ centered at $x_{0}$,
\item In the annular region $h^{\delta_{0}} \leq |x-x_{0}| \leq R_{0}$, take a covering by balls of size $h^{\delta_{0}+\delta_{1}}$.
\end{itemize}  
For $\chi_{p}$ supported in this annular region, $\kappa(p)$ will behave as $h^{-\delta_{0}}$, whereas $\eps(p)$ will be of size $h^{\delta_{0}+\delta_{1}}$. This will be enough to have a small error term. 

In dimension $n\geq3$, the key is to iterate this procedure according to the stratification of the corner domain, and we get a suitable covering of $\overline{\Omega}$, as follows (see \cite[Lemma]{BruPof16} and \cite[Appendix B]{BoDauPof14}): $\overline{\Omega}\subset \cup_{p\in \cP}\cB(p,h^{\delta_{0}+\ldots+\delta_{k}})$, where the ball $\cB(p,2h^{\delta_{0}+\ldots+\delta_{k}})$ is contained in a map neighborhood of $p$, and the curvature associated with this map-neighborhood satisfies 
\bel{E:bdcurv}
\kappa(p) \leq \frac{c(\Omega)}{h^{\delta_{0}+\ldots+\delta_{k-1}}}.
\ee
In this procedure, $k$ depends on the point $p$. Moreover, $k$ can be chosen between 0 and $\nu$, where $\nu$ is the smallest integer satisfying
$$\forall \dx \in \gC(\Omega), \quad l(\dx) \geq \nu \implies \Pi_{\dx}\ \  \mbox{is polyhedral}.$$
Note that $\nu \leq n-1$, and that $\nu=1$ if and only if $\Omega$ is polyhedral. 

We will take a partition of unity $(\chi_{p})_{p\in \cP}$ satisfying $\supp(\chi_{p})\subset B(p,2h^{\delta_{0}+\ldots+\delta_{k}})$, and
\beq
\label{E:uniformestimatesgradient}
\left\{
\begin{aligned}
& \sum_{p\in \cP} \chi_{p}^2=1  \quad \mbox{on} \ \ \overline{\Omega},
\\
& \sum_{p\in \cP} \|\nabla \chi_{p}\|^2_{\infty} \leq C(\Omega)h^{-2\delta} \ \ \mbox{with} \ \ \delta=\delta_{0}+\ldots+\delta_{\nu}.
\end{aligned}
\right.
\ee
We deduce from \eqref{IMSrem}
\beq
\label{E:IMSschema}
L \geq h \sE(\Omega)+h\sum_{p\in \cP}R(\eps(p),\kappa(p))+\cO(h^{2-2\delta}).
\ee
This preliminary step will provide a lower bound if 
\begin{itemize}
\item The lowest local energy is reached, in order to avoid degeneracy such as $\sE(\Omega)\in \{0,-\infty\}$. This will be guaranteed by the application of Theorem \ref{T:Fsci} to the local energy. 
\item The scales are chosen such that $\delta \in(0,\frac{1}{2})$ and $R=o(1)$, for this last estimate, we will use \eqref{E:bdcurv}. 
\end{itemize}

\chapter{The Robin Laplacian}
\label{C:robin}
In this chapter, we consider the Laplacian in a bounded corner domain $\Omega\in \gD(M)$, where $M$ is a Riemannian manifold (mainly, $M=\R^{n}$ or $M=\dS^{n}$), with a Robin boundary condition $\partial_{\nu}u-\alpha u=0$, and investigate the asymptotics of its first eigenvalues as $\alpha\to+\infty$. This is equivalent to the semi-classical regime $h^{1/2}\partial_{\nu}-u=0$, as $h\to0$. In the Section \ref{S:robinlite}, we present the operator and give a short overview of the asymptotics of its first eigenvalue from the literature. In Section \ref{SS:Robincoin}, using the tools developed in Chapter \ref{C:corner}, we present the recursive analysis leading the asymptotic behavior of the first eigenvalue in a corner domain, including an {\it a priori } two-side estimate, and the analysis of the essential spectrum of the tangent operators. We also review our results in the case where the Robin boundary condition have a vanishing Dirichlet weight function, leading to a non self-adjoint operator. In Section \ref{SS:RobinReg}, we give a more precise asymptotics of the low-lying eigenvalues in a regular domain, using an effective operator, defined on the boundary, and involving the mean curvature. We apply our result to show a reverse Faber-Krahn inequality. 

\section{The Robin Laplacian with large Dirichlet parameter}
\label{S:robinlite}


We consider the Robin eigenvalue problem on a bounded corner domain $\Omega$: 
\beq
\label{E:evproblem}
\left \{ 
\begin{aligned}
&-\Delta u =\lambda u \ \ \mbox{on}\ \  \Omega,
\\
&\partial_{\nu}u-\alpha u=0 \ \ \mbox{on}\ \  \partial\Omega,
\end{aligned}
\right. 
\ee
where $\alpha$ is a real parameter. 
We denote by $\cQ_{\alpha}^{\Omega}$ the associated quadratic form:
\beq
\label{E:DefcqQ}
\cQ_{\alpha}^{\Omega}(u):=\| \nabla u\|_{L^2(\Omega)}^2-\alpha \|u\|_{L^2(\partial\Omega)}^2, \quad u\in H^1(\Omega).
\ee
Since $\Omega$ is bounded and is the finite union of Lipschitz domains (see \cite[Lemma A.A.9]{Dau88}), the trace map from $H^{1}(\Omega)$ into $L^{2}(\partial\Omega)$ is compact and the quadratic form $\cQ_{\alpha}^{\Omega}$ is lower semi-bounded. We define its the self-adjoint operator $\cL_{\alpha}^{\Omega}$ as its Friedrichs extension, whose spectrum is a sequence of eigenvalues $(\lambda_{j}^{R}(\alpha,\Omega))_{j\geq1}$ (shorted to $\lambda_{j}^{R}(\alpha)$ when there is no ambiguity on the domain $\Omega$), in particular $\lambda_{1}^{R}(\alpha)$ is called the principal eigenvalue of the system \eqref{E:evproblem}. 

It is well known that $\lambda_{1}^{R}(\cdot)$ is decreasing, concave, that $\lambda_{1}^{R}(0)=0$ and that $\lim_{\alpha\to-\infty}\lambda_{1}^{R}(\alpha)$ is the first Dirichlet eigenvalue on $\Omega$, but the limit as $\alpha\to+\infty$ appears to be more singular, in the sense that it is not clear what limit problem would drive the asymptotics, and therefore, do not enter the framework of regular perturbation theory.

Using a constant test functions, it is direct to see from the min-max principle that $\lambda_{1}^{R}(\alpha)\to-\infty$ as $\alpha\to+\infty$. Clearly, the limit as $\alpha\to+\infty$ is linked to the limit $h\to0^+$ for the Laplacian with the boundary condition $h^{1/2}\partial_{\nu}u-u=0$. This asymptotics regime for the Robin Laplacian has several application in reaction diffusion systems (\cite{LaOckSa98}), surface superconductivity (\cite{MonIn00,GiSm07,AsBaPP15}), and the study of its spectrum has receive a lot of interests since then (\cite{Dan06,LevPar08,ExMinPar13,Pank13,CakChauHad,FreKr15,Dan13,HePank15,HeKaRay17,KhalPank16,Kha18}).

Therefore, we are able to define our operator in this chapter: 
\begin{notation}
\label{N:robin}
In this Chapter, the operator $L_{h}$ is defined by $h^{-1}L_{h}=\alpha^{-2}\cL^{\Omega}_{\alpha}$, with $\alpha=h^{-1/2}$. In particular, the quadratic for associated with $L_{h}$ is 
$$Q_{h}:u\mapsto h^2\int_{\Omega} |\nabla u|^2-h^{3/2}\int_{\partial\Omega}|u|^2.$$
Therefore, $\lambda_{j}(h)=\alpha^{-4}\lambda_{j}^{R}(\alpha)$. For a cone $\Pi\in \gP_{n}$, the associated tangent operator $L_{h}[\Pi]$ is defined as the extension of the form.
$$
\cQ_{h}[\Pi]:u \mapsto h^2\| \nabla u\|_{L^2(\Pi)}^2-h^{3/2}\|u\|_{L^2(\partial\Pi)}^2, \quad u\in H^1(\Pi), $$
and the normalized model operator is $L[\Pi]=L_{1}[\Pi]$. The local energy $\gC(\Omega)\ni \dx \mapsto E(\Pi_{\dx})\in \R$ is now well defined by \eqref{E:DE}, as the bottom of the spectrum of $L[\Pi_{\dx}]$.
\end{notation}
Of course, to be totally rigorous, we have to show that the quadratic form, involved in the above definition of $L[\Pi]$, is bounded from below. This is a non-trivial fact which is linked to the recursive definition of $\gP_{n}$, see Theorem \ref{T:essentialRobincone}.

Let us notice that the text presents two equivalent, slightly confusing, notations: one with the parameter $h\ll1$, and another one with the parameter $\alpha \gg1$. We did so because the first one matches with the general framework presented in the other chapters, and the second one is the most usual one in the literature. In this chapter, we will present our asymptotics for the quantity $\lambda_{j}^{R}(\alpha)$.

Using the scaling $x\mapsto h^{-1/2} x$ in a cone for the Rayleght quotient $\frac{Q_{h}(u)}{\|u\|^2}$, we see that we are within the framework of Section \ref{S:cornerpure} since \eqref{E:homoge} is satisfied. In a bounded domain $\Omega\in \gD(\R^{n})$, if we accept the fact that the eigenfunctions tend to be localized near some part of the boundary, according to the above scaling, as $\alpha\to+\infty$, then it is expected that 
 \beq\label{E:LevPar}
 \lambda_{1}^{R}(\alpha) \underset{\alpha\to+\infty}{\sim} C(\Omega)\alpha^2.
 \ee
Let us make a short overview of the pre-existing results. For smooth domains, this is proved with $C(\Omega)=-1=E(\R^{n}_{+})$ (see \cite{LaOckSa98,LouZhu04} and \cite{DaKe10} for higher eigenvalues), with various improvements depending on the geometry of the boundary, see Section \ref{SS:RobinReg}. 

The sectors of opening $\theta\in (0,2\pi)$, denoted by $S_{\theta}$, enjoy an explicit expression for their ground state energy: 
\beq
\label{E:Esectors}
 E(S_{\theta})=\left\{ \begin{aligned} &-\sin^{-2}\tfrac{\theta}{2} \ \ \mbox{if} \ \ \theta\in (0,\pi) \, ,\\ &-1 \ \ \mbox{if} \ \ \theta\in (\pi,2\pi)\, . \end{aligned}\right.
 \ee
For planar polygonal domains with corners of opening $(\theta_{k})_{k=1,\ldots,N}$, it is conjectured in \cite{LaOckSa98} that \eqref{E:LevPar} holds with
$$C(\Omega)= - \max_{ 0<\theta_{k}<\pi}(1,\sin^{-2} \tfrac{\theta_{k}}{2})=\min_{\theta_{k}} E(S_{\theta_{k}}).$$
Therefore, the asymptotics \eqref{E:LevPar} seems to hold, at least for regular domains and two-dimensional polygons, with 
\beq
\label{E:Comega}
C(\Omega)=\inf_{x\in \partial\Omega} E(\Pi_{x})=\sE(\Omega),
\ee
which is the same formulation that Conjecture \eqref{E:megameta}.  This is finally proved in \cite{LevPar08} for domains with corners satisfying the uniform interior cone condition, although the finiteness of $\sE(\Omega)$ is not studied.

\section[Asymptotics for the first eigenvalue]{Asymptotics for the first eigenvalue (including the $\delta$-interaction on hypersurfaces)}
\label{SS:Robincoin}
As suggested in the introduction, see Section \ref{S:intro1}, the asymptotics \eqref{E:LevPar}, with \eqref{E:Comega}, rises two questions: Is $\sE(\Omega)$ finite, and is it possible to have an a priori remainder? These questions are in fact quite related, and the first one will get a positive answer if the infimum in \eqref{E:Comega} is reached. Our aim is therefore to apply Theorem \ref{T:Fsci}. Next, our goal is to prove {\it a priori} remainder estimates. Finally, we have proved the following theorem, whose proof, from \cite{BruPof16}, is inspired by the strategy developed in \cite{BoDauPof14}: 

\begin{theorem}
\label{T:Robin}
Let $\Omega \in \gD(M)$, and $\sE(\Omega)$ be defined in \eqref{D:lle}. Then 
\begin{enumerate}
\item
$\sE(\Omega)>-\infty$, and there exists $x_{0}\in \partial\Omega$ such that $\sE(\Omega)=E(\Pi_{x_{0}})$.
\item There exist $\alpha_{0}\in \R$, two constants $C^{\pm}>0$, and two integers $0\leq \nu' \leq \nu\leq n-2$, such that 
\bel{E:metarobin}
\forall \alpha\geq \alpha_{0}, \quad -C^{-}\alpha^{2-\frac{2}{2\nu+3}}\leq \lambda_{1}^{R}(\alpha)-\alpha^2\sE(\Omega) \leq C^{+}\alpha^{2-\frac{2}{2\nu'+3}}.
\ee
\end{enumerate}
\end{theorem}
We now give the main lines of the proof of the Theorem:
\paragraph{Lower bound}  
Our operator enters the framework described in Section \ref{S:Lowerbound}, indeed using $h=\alpha^{-2}$, the quadratic form $Q_{\alpha}^{\Omega}$ can be written 
$$H^{1}(\Omega)\ni u\mapsto h^{-2}\left(\int_{\Omega}h^2|\nabla u|^2-h^{3/2}\int_{\partial\Omega}|u|^2\right).$$
The error term coming from the approximation of the metrics in \eqref{E:IMS1} satisfies $R=\cO(\eps(p)\kappa(p))=h^{\delta_{k}}$, see \cite[Section 3 \& 4]{BruPof16}. Therefore \eqref{E:IMSschema} provides 
$$L_{h} \geq h \sE(\Omega)+\sum_{k=0}^{\nu}\cO(h^{1+\delta_{k}})+\cO(h^{2-2\delta}).$$
The integer $\nu$ depends on $\Omega$ as follows: 
$$\nu=\inf\{m\in \N, \ \mbox{all chains of length }m\mbox{ are polyhedral} \}.$$
In particular, $\nu=0$ if and only if $\Omega$ is polyhedral. The optimization of remainders is done by choosing $\delta_{0}=\ldots=\delta_{\nu}$ and $2-\delta_{0}=2\delta=2(\nu+1)\delta_{0}$, that is $\delta_{0}=\frac{2}{2\nu+3}$. Linking $L_{h}$ and $\cL_{\alpha}^{\Omega}$, we deduce from the min-max principle that there exist $\alpha_{0}\in \R$ and $C^{-}>0$ such that
 \beq
 \label{E:lowerboundproved}
 \forall \alpha \geq \alpha_{0}, \quad \lambda_{1}^{R}(\alpha) \geq \alpha^2 \sE(\Omega)-C^{-}\alpha^{2-\frac{2}{2\nu+3}}.
 \ee

\paragraph{Bottom of the spectrum of the tangent operator}
Let $\Pi \in \gP_{m}$, and let $\Gamma$ be its reduced cone. In some suitable coordinates, we may write
\bel{E:PIdecompose}
\Pi=\R^{m-n}\times \Gamma
\ee
with $\Gamma\in \gP_{n}$ an irreducible cone and $n \leq m$. The associated Robin Laplacian admits the following decomposition: 
\beq
\label{E:decomtensor}
L[\Pi]=-\Delta_{\R^{m-n}}\otimes \Id_{n}+\Id_{m-n}\otimes L[\Gamma]
\ee
In particular 
$$\spec(L[\Pi])=[E(\Gamma),+\infty).$$
We denote by $\omega$ the intersection of $\Gamma$ with the unit sphere. Then we have the intermediate result, stated in \cite{BruPof16}:
\begin{theorem}
\label{T:essentialRobincone}
Assume that $\sE(\omega)>-\infty$. Then 
\begin{enumerate}
\item $E(\Pi)>-\infty$, and the Robin Laplacian, $L[\Pi]$, is well defined as the Friedrichs extension of $\cQ[\Pi]$, with form domain $H^{1}(\Pi)$,
\item Assume moreover that $\Pi$ is irreducible. Then the bottom of the essential spectrum of $L[\Pi]$ is $\sE(\omega)$.
\end{enumerate}
\end{theorem}
This theorem relies on the fact that 
\beq
\label{E:fqpolaire}
\cQ[\Pi](u)=\int_{r>0} \left(|\partial_{r}u|^2+\frac{1}{r^2}\cQ_{r}^{\omega}(u(r,\cdot))\right) r^{n-1} \rd r,
\ee
where $(r,\theta)\in (0,+\infty) \times \omega$ denotes spherical coordinates in the cone $\Pi$. We now see that in $\Pi$, for each $R$ fixed, the Robin Laplacian in the region $r \approx R$ is linked to the Robin Laplacian with parameter equal to $R$. As $R$ gets large, we can combine \eqref{E:fqpolaire} and the hypothesis that $\sE(\omega)$ is finite in order to have a lower bound for the quadratic form. The Persson's Lemma provides the bottom of the essential spectrum. Let us notice that the the quantity giving the bottom of the essential spectrum is an infimum over singular chains, and that its structure reminds of the HVZ's Theorem for the $N$-body problem, see \cite{GeoNis17} for a recent approach.

Note that, at first glance, the second point is of interest only when $\Pi$ is irreducible. However, if $\Pi$ writes as in \eqref{E:PIdecompose}, and if $\omega$ is the section of $\Gamma$, the quantity $\sE(\omega)$ will play the role of a threshold, as we will show below.

\paragraph{Finiteness of the lowest local energy}
For the Robin Laplacian, this is done by induction on the dimension $n$. We assume that $\sE(\omega)$ is finite for all $\omega\in \gD(M)$, where $M$ is an $n-1$ dimensional manifold without boundary. We consider a $n$ dimensional corner domain $\Omega$. First, by standard perturbation argument, the bottom of the spectrum is regular with respect to the geometry, in the sense that $\gP_{n}\ni \Pi\mapsto E(\Pi)$ is continuous with respect to the distance $\dD$ introduced in Definition \ref{D:Ordrechain}.

We now prove the monotonicity of the local energy. 
Let $x_{0}\in \overline{\Omega}$. Denote by $\omega_{x_{0}}$ the section of the reduce cone of $\Pi_{x_{0}}$, then, we need to introduce {\it the second energy level}, 
\bel{D:Estar0}
E^{*}(\Pi_{x_{0}})=\inf_{x_{1}\in \partial\omega_{x_{0}}}E(\Pi_{x_{0},x_{1}}),
\ee
where the tangent structure $\Pi_{x_{0},x_{1}}$ associated with a chain $(x_{0},x_{1})$ has been defined in Section \ref{S:cornerpure}. By definition of $\sE(\cdot)$, 
\bel{D:Estar}
\sE(\omega_{x_{0}})=E^{*}(\Pi_{x_{0}}).
\ee
Therefore, since $\omega_{x_{0}}\in \gD(\dS^{n-1})$, the recursive hypothesis combined with  Theorem \ref{T:essentialRobincone} show that
$$\forall x_{1}\in \overline{\omega_{x_{0}}}, \quad -\infty< E(\Pi_{x_{0}}) \leq E(\Pi_{x_{0},x_{1}}),$$
We deduce by immediate recursion that $\gC(\Omega)\ni\dx\mapsto E(\Pi_{\dx})$ is monotonous. As a consequence of Theorem \ref{T:Fsci}, $\partial\Omega\ni x \mapsto E(\Pi_{x})$ is lower semi-continuous, therefore it reaches its infimum, and it is finite. This concludes the induction, and proves the first point of Theorem \ref{T:Robin}.

Given an irreducible cone $\Pi$, one may ask for a condition such that the inequality $E(\Pi)\leq E^{*}(\Pi)$ is strict, i.e. wether there exists discrete eigenvalue below the essential spectrum of $L[\Pi]$. This question, interesting in itself, receives an answer in the following cases: 
\begin{itemize}
\item If $\Pi$ is a sector of opening $\alpha$, the answer is yes if and only if $\alpha\in (0,\pi)$, i.e. if and only if it is a convex sector, see \eqref{E:Esectors}.
Moreover, it is proven in \cite{KhalPank16} that when $\alpha\in (0,\pi)$, there exist only a finite number of eigenvalues below the essential spectrum. 
\item  If the complementary of the cone is convex, then the answer is no (\cite[Corollary 3]{Pank16}).
\item If the mean curvature is positive at one point of the boundary of $\omega:=\Pi\cap \dS^{n-1}$, then the answer is yes. Moreover, there exists an infinite number of eigenvalues below the essential spectrum (\cite[Theorem 6]{Pank16}).
\item In \cite{BruPankPof16}, we assume that $n=3$ and that $\omega$ is smooth. Denote by $\kappa$ the geodesic curvature of $\partial\omega$, and by $\kappa_{+}$ its positive part. Then $\sE(\Pi)=-1$, and the number of eigenvalue of $L[\Pi]$ in $(-\infty,-1-\lambda)$ behaves, as $\lambda\to0^{+}$, as 
$$\frac{1}{8\pi\lambda}\int_{\partial\omega} \kappa_{+}(k) \rd s.$$
\item When $\omega$ is convex, upper and lower bound on $E(\Pi)$ are given in \cite[Section 5]{LevPar08}, using geometrical quantities. These bounds become an exact value when $\Pi$ is polygonal and $\omega$ admits an inscribed circle.
\end{itemize}

\paragraph{Construction of test-functions and upper bound}
We aim to construct a test-functions for the operator in $\Omega$, localized near a minimizer $x_{0}\in \partial\Omega$ of the local energy. As described in the last paragraph, if $E(\Pi_{x_{0}})<E^{*}(\Pi_{x_{0}})$, then there exists an eigenfunction (with exponential decay) $\psi_{0}$ for $L[\Gamma_{x_{0}}]$. Next, we consider the function $\one \otimes \psi_{0}\in L^{\infty}(\R^{n-d}\times \Gamma_{x_{0}})$, where $d$ is the reduced dimension of $\Pi_{x}$. A quasi-mode for $L[\Pi_{x_{0}}]$ is easily constructed from this function, after suitable rotation, cut-off and scaling. This quasi-mode is then transported by a local map in a quasi-mode in $\Omega$, and estimate of the remainder provide an upper bound of the first eigenvalue. It is called a {\it sitting quasi-mode}, in the terminology of \cite{BoDauPof14}. This procedure is rather standard.

But if $E(\Pi_{x_{0}})=E^{*}(\Pi_{x_{0}})$, the operator on the reduced cone $\Gamma_{x_{0}}$ does not have discrete eigenvalues below its essential spectrum, and it is not clear how to localize a test function with energy $E(\Pi_{x_{0}})$. The method is to construct a quasi-mode qualified as {\it sliding} in \cite{BoDauPof14}. Notice that the spirit close to the construction of quasi-modes in manifold with corners, described in \cite{Gri17}. We describe the main lines of the proof below.

First, there always exists a chain $\dx=(x_{0},\ldots,x_{\nu})$ initiated at $x_{0}$ such that $E(\Pi_{\dx})=E(\Pi_{x_{0}})=\sE(\Omega)$, and $E(\Pi_{\dx})<E^{*}(\Pi_{\dx})$. This is shown by recursion, see \cite[Proposition 7.1]{BruPof16}, starting with the half-plane, where $E^{*}=0$ and $E=-1$. Therefore it is possible to construct a test-function for $L[\Pi_{\dx}]$. 

Next, this function is then used through a sequence of scaling and translations in order to get a test function defined in $\Pi_{x_{0}}$, localized a neighborhood of $0$, but whose support avoids the origin (and it also avoids all the points $x_{k}$, $k\leq \nu-1$, in the recursive steps). The sizes of the scaling are different, and not determined at this stage.


The last step is similar to the sitting case: we transport the test function, defined in $\Pi_{x_{0}}$, into a function in $\Omega$. Roughly, the energy of this last test function is $\alpha^2 E(\Pi_{\dx})$, modulo some remainders, coming from the use of cut-off and approximations of the metric. The scales are then set in order to optimize the remainders. The min-max-principle provides the upper bound of \eqref{E:metarobin}.

\paragraph{The Laplacian with a strong $\delta$-interaction}
Our results are true for the Laplacian acting on a hypersurface with corners with a strong $\delta$-interaction. 
 Let $\Omega \in \gD(M)$ be a corner domain and let $S=\partial \Omega$ be its boundary. We consider $\cL^{S,\delta}_{\alpha}$ the self-adjoint extension associated with the quadratic form
$$u\mapsto \| \nabla u\|_{L^2(M)}^2-\alpha \|u\|_{L^2(S)}^2, \quad u\in H^1(M).$$
The associated boundary problem is the Laplacian with the derivative jump condition  across the closed hypersurface $S$: $[\partial_{\nu} u]_{\partial \Omega} = \alpha u$.
When $M=\R^{n}$, it is well known (see e.g. \cite{BrExKuSe94}) that since $S$ is bounded, $\cL^{S,\delta}_{\alpha}$ is a relatively compact perturbation of $-\Delta$ on $L^2(\R^n)$ and then
$$\sigma_{ess}( \cL^{S,\delta}_{\alpha}) = [0, + \infty).$$
Moreover $\cL^{S,\delta}_{\alpha}$ has a finite number of negative eigenvalues. If we denote by $\lambda^{\delta}_{1}(\alpha)$ the lowest one, by applying the strategy developed for the Robin Laplacian, all the above results are still valid replacing $\lambda^{R}_{1}(\alpha)$ by $\lambda_{1}^{\delta}(\alpha)$.
We still define the tangent operator at a point $x$ as $\cL^{S_{x},\delta}_{1}$, where $S_{x}$ is the boundary of $\Pi_{x}$. The associated local energy at $x$, $E^{\delta}(S_{x})$, is the bottom of the spectrum of $L^{S_{x},\delta}_{1}$, and their infimum is $\sE^{\delta}(S)$. 
Then, as explained in \cite{BruPof16}:
\begin{theorem}
Theorems \ref{T:Robin} and \ref{T:essentialRobincone} remain valid when replacing the Robin Laplacian $\cL^{\Omega}_{\alpha}$ by the $\delta$-interaction Laplacian $\cL^{S,\delta}_{\alpha}$. 
\end{theorem}
When $x$ belongs to the regular part of $S$, $S_{x}$ is an hyperplane and 
\beq
\label{E:}
E^{\delta}(S_{x})=-\tfrac14,
\ee
see \cite{ExYo02}. Therefore $\sE^{\delta}(S)=-\frac14$ when $S$ is regular, and we recover the known main term of the asymptotic expansion of $\lambda^{\delta}_{1}(\alpha)$ proved in dimension $2$ or $3$ (see \cite{ExYo02,  ExPa14, DiExKuPa15}). 

To our best knowledge the only studies for $\delta$-interactions supported on non smooth hypersurfaces are for broken lines and conical domains with circular section (see \cite{BerExLot14, DuRa14,ExKon15,LotOur16}). In that case, it is proved in the above references that the bottom of the essential spectrum of $L^{S,\delta}_{\alpha}$ is $- \alpha^2/4$, which can be deduced from our Theorem \ref{T:essentialRobincone} together with the scaling argument. In view of our result, this remains true  when the section of the conical surface is smooth, and we are able to compute the bottom of the essential spectrum for a wider class of $\delta$-interaction on cones.

Moreover, our work seems to be the first result giving the main asymptotic behavior of $ \lambda^{S,\delta}_{\alpha}$ for interactions supported by general closed hypersurfaces with corners.

\begin{remark}
\label{R:weight}
For the Robin Laplacian and  the $\delta$-interaction Laplacian, we can add a smooth positive weight function $G$ in the boundary conditions. These conditions become, for the Robin condition $\partial_\nu u = \alpha G(x)  u$ and for the $\delta$-interaction case, $[\partial_{\nu} u ] = \alpha G(x)  u$. In our analysis, for $x \in \partial \Omega$ fixed, we change $\alpha$ into $\alpha G(x)$ and clearly, the results are still true by replacing  $\sE(\Omega)$ and $\sE^\delta(S)$ by:
$$ \sE_G( \Omega):=\inf_{x\in \partial \Omega}G(x)^2 E(\Pi_{x}), \qquad \sE_G^\delta(S):=\inf_{x\in S}G(x)^2 E^\delta (S_{x}).$$
For the Robin Laplacian, such cases were already considered in \cite{LevPar08} and \cite{CoGa11}.
\end{remark}
Note that the Laplacian with a Robin-type boundary condition $\omega\partial_\nu u =  u$, involving a variable function $\omega$ which vanishes at a point $x_{0}$, can be very different. In dimension 2, denote by $s$ an arclength parameter of $\partial\Omega$ near $x_{0}$, and assume that $\omega$ vanishes at order 1 (i.e. $\omega (s_{0})=0$ and $\omega'(s_{0})\neq0$). The associated quadratic form $u\mapsto \|\nabla u\|^2_{L^2(\Omega)}-\||\omega|^{-1/2}u\|_{L^{2}(\partial\Omega)}$ is not bounded from below in $H^{1}$, and it is not clear how to define a self-adjoint operator. This problem has been noticed in the model case of a half-disc, \cite{BerDen08,MarRoz09}. Using the Kondratiev theory (\cite{Kon67}), we are able to show in \cite{NazPof18} that the Robin Laplacian has indices of defect equal to $(1,1)$ (\cite{ReSi78}). Therefore, the spectrum of its adjoint covers the whole complex plane $\C$, and we are able to describe its self-adjoint extensions as a one-parameter family $L(\theta)$, for $\theta\in \dS^{1}$. We also study the dependency of their eigenvalues with respect to $\theta$. The situation where $\omega$ vanishes to other orders may be very different, and may not be covered by the theory of Kondratiev. We plan to investigate these cases.

\section{Effective Hamiltonian in the regular case: the role of the curvature}
\label{SS:RobinReg}
We now assume that $\Omega$ is a $C^{2}$ domain, and we describe how to get a more precise asymptotics of the low-lying eigenvalues. As the analysis of the last section has shown, the minimizer of the local energy will govern the first term of the asymptotics. But for a regular domain, $E(\Pi_{x})=-1$ for all $x\in \partial\Omega$. If one thinks of the local energy as an effective potential, then it would say in that case that its extremum is reached on the whole boundary; in this sense the problem is not very far of the {\it puits sous vari\'et\'es} for the harmonic approximation described in \cite{HeSj87}, and it is expected that an effective operator, defined on $\partial\Omega$, will lead the asymptotics of the low-lying eigenvalues of $\cL_{\alpha}^{\Omega}$, as $\alpha\to+\infty$.

In the literature, first results are obtained in dimension $n=2$: if we denote by $\gamma:\partial\Omega \to \R$ the curvature of the boundary, and $\gamma_{\max}$ its maximum, then it is proved in \cite{ExMinPar13,Pank13} that, as $\alpha\to+\infty$, for fixed $j\in \N$,
\beq
\label{A:Robin2d}
  \lambda_{j}^{R}(\alpha)=-\alpha^2-\gamma_{\max}\alpha+\cO(\alpha^{2/3}) .
\ee
Note that in case of balls and spherical shells, these asymptotics are done through analytical ODEs, see \cite{FreKr15}.

The asymptotics \eqref{A:Robin2d} rises several questions: what is the analogous in higher dimension, and is it possible to see the influence of $j\in \N$ in the next term of the asymptotics? In \cite{HeKa15}, the authors assume that $\Omega\subset \R^{2}$ is $C^{\infty}$, and that $\gamma$ admits a unique non-degenerate maximum. They prove that this maximum acts as a wells and they provide a full asymptotic expansion of the eigenvalues (in the form \eqref{curvaturewell} below)).

For a regular domain $\Omega\subset \R^{n}$, we denote by $H:\partial\Omega\to$ the mean curvature of its boundary and by $-\Delta_{S}$ the (positive) Laplace-Beltrami operator on the hypersurface $\partial\Omega$. Then, we prove in \cite{PankPof16}: 
\begin{theorem}
Assume that $\Omega$ is $C^{3}$. Denote by $L^{\eff}_{\alpha}$ the operator $-\Delta_{S}-(n-1)\alpha H$, acting on $H^{2}(\partial\Omega)$, and $\mu_{j}(\alpha)$ its $j$-th eigenvalue. 
Let $j\in \N$ be fixed, then, as $\alpha\to+\infty$: 
\beq
\label{A:robingene}
 \lambda_{j}^{R} (\alpha)=-\alpha^2+\mu_{j}(\alpha)+O(1).
 \ee
\end{theorem}
As $\alpha\to+\infty$, the operator $L^{\eff}_{\alpha}$ has the form of a semiclassical Schr\"odinger operator on $\partial\Omega$, the potential being proportional to $-H$. In particular its eigenvalues satisfy, as $\alpha$ gets large,  $\mu_{j}(\alpha)\sim -(n-1)\alpha H_{\max}$, where $H_{\max}$ denotes the maximum of $H$. More precise asymptotics enters the framework of harmonic approximation, see \cite{HeSj84,Si82,DiSj99}: hypotheses on $H$ near its minimum will imply more structure in the asymptotics. In particular: 
\begin{corollary}
\label{C:asymptRobfirsterm}
There holds, as $\alpha\to+\infty$: 
$$\lambda_{j}^{R}(\alpha)=-\alpha^2-(n-1)\alpha H_{\max}+o(\alpha)$$
 Assume moreover that $\Omega$ is $C^{\infty}$, and that $H$ admits a unique global maximum at $s_{0}$ and that the Hessian of $-(n-1)H$
at this maximum is positive-definite. Denote by $\mu_{k}$ its eigenvalues
and set
\[
\cV=\bigg\{ \sum_{k=1}^{\nu-1}\sqrt{\frac{\mu_{k}}{2}}\,\big(2n_{k}-1\big), n_{k}\in \N\bigg\}.
\]
Let us sort the elements of $\cV$ in the increasing order, repeating the terms if they appear multiple
times, and denote by $e_{j}$ the $j$-th element. Then for each $j\in\N$ there holds, as $\alpha\to+\infty$,
\bel{curvaturewell}
\lambda_{j}^{R}(\alpha)=-\alpha^2-(n-1)\alpha H_{\max}+e_{j}\alpha^{1/2}+\cO(\alpha^{1/4}).
\ee
Moreover, if $e_{j}$ is of multiplicity one, the remainder estimate can be improved to $\cO(1)$. 
\end{corollary}
Notice that various improvements of the above result exist, depending on some changes in the assumptions: 
\begin{itemize}
\item An analogous of the result holds when $\Omega$ is not compact, provided it has a good structure at infinity, see \cite[Definition 1.1]{PankPof16}, in the spirit of \cite{CarExKre04}.
\item In dimension 2, if $H$ has a unique maximum of order $2<p<+\infty$, a three-terms asymptotics is still available, see \cite[Corollary 1.8]{PankPof16}.
\item If the domain $\Omega$ is only $C^{2}$, \eqref{A:robingene} still holds with a remainder in $\cO(\log \alpha)$ instead of $\cO(1)$.
\item In dimension $n=2$, assuming that there exist a unique non-degenerate curvature well, \eqref{curvaturewell} has been obtained in \cite{HeKa15}, with a full asymptotic expansion in power of $\alpha^{1/2}$. In the case of a symmetric domain with two curvature wells, the tunneling effect is analyzed in \cite{HeKaRay17}, in particular it is proved that the tunneling effect is given by the one of $L^{\eff}_{\alpha}$, which can be seen as a Schr\"odinger operator acting on $\dS^{1}$ with a double-wells potential.
\end{itemize}

Following this analysis, the question of a refinement of \eqref{E:metarobin} (and of an asymptotics for the higher eigenvalues) for non-smooth boundary is a relevant question. For polygons, this is the object of \cite{Khalile}. Following the ideas of \cite{Bon06,BonDau06}, this problem is treated in \cite{Kha18}: if each model problem at a vertex $v_{j}$ has $K_{j}$ eigenvalues below its essential spectrum, then the first $K:=\sum_{j}K_{j}$ eigenvalues will be ``attracted'' by the corners, with an exponentially small interaction. The asymptotics expansion of the next eigenvalues requieres a more global approach, since the sides of the polygons will now contribute. This is done in \cite{KhaOurPank18} (see also \cite{Pank15}), under an additional hypothesis on the essential spectrum of the model problems: in a polygon with zero curvature, once the corners have ``attracted'' the $K$ first eigenvalues, the $K+j$ eigenvalue of the Robin Laplacian satisfies
$$\lambda_{K+j}(\alpha)=-\alpha^2+\mu_{n}+\cO\big(\frac{\log \alpha}{\alpha}\big)$$
where $\mu_{n}$ is the $n$-th eigenvalue of the Laplacian on the graph formed by the sides of the boundary, with a Dirichlet boundary condition at each junction. Note that this approach is inspired by \cite{Gri08} for a similar problem.

\paragraph{Application to a Faber-Krahn inequality}
The optimization of eigenvalues under geometrical constraints has received a lot of interest since more than a century, originating from Lord Rayleigh's {\it Theory of sound}. He conjectures the following property, now classical: : among bounded domains $\Omega$ of fixed volume, the first eigenvalue of the Dirichlet Laplacian should be minimized by the ball. This has been proven by Faber and Krahn in the 1920's. This kind of questions has been extended to numerous optimization problems, such as the second Neumann eigenvalue and higher Dirichlet eigenvalues, forming the family of  {\it isoperimetric spectral inequalities}. We refer to \cite[Sections 3-7]{Henrot} for an overview containing the historical references. 

In this part, we note $\lambda_{1}^{R}(\alpha,\Omega)$ the first eigenvalue of the Robin Laplacian, in order to emphasize the dependency on the domain. It has been proved in \cite{Boss86} in dimension 2, and in \cite{Dan06} in any dimension, that $\lambda_{1}^{R}(\alpha,\Omega)$ is also minimized by a ball when $\alpha$ is fixed in $(-\infty,0)$.  

But, according to the value of the derivative $\partial_{\alpha}\lambda_{1}^{R}$ when $\alpha=0$, it has been conjectured in \cite{Bar77} that the ball should become the maximizer for $\lambda_{1}^{R}(\alpha,\Omega)$ when $\alpha>0$ is fixed. 
This {\it reverse} Faber-Krahn inequality has been proved in dimension $n=2$ for $\alpha$ small enough (\cite[Section 4]{FreKr15}), but, surprisingly, this conjecture has been disprove for large $\alpha$: in \cite[Section 3]{FreKr15}, it is proved that if $\cB$ and $\cA$ are a ball and a spherical shell of same volume, then, for $\alpha$ large enough, $\lambda_{1}^{R}(\alpha,\cB)<\lambda_{1}^{R}(\alpha,\cA)$. Such a counter-example relies on explicit calculations, using the spherical invariance of the domains in order to express the eigenvalues as roots of special functions. 

In view of Corollary \ref{C:asymptRobfirsterm}, the maximization of $\lambda_{1}^{R}(\alpha,\Omega)$ leads to the following question: 
\beq
\label{Q:meancurv}
\mbox{Set the volume of } \Omega, \mbox{how to minimize } H_{\max} \, ? 
\ee
The counter example to the Faber-Krahn inequality shows also that this problem has no solution without additional constraint, indeed a thin annulus of large radius can have a fixed volume, but its mean curvature can be arbitrary small. In \cite{PankPof15}, we prove:

\begin{theorem}
\label{T:optimMean}
Let $\Omega\subset\R^{n}$ be a bounded star-shaped regular domain, and $B$ a ball of same volume. Then   
$$H_{\max}(\Omega) \geq H_{\max}(B),$$
moreover there is equality if and only if $\Omega$ is a ball.
\end{theorem}
The above theorem relies on the standard isoperimetric inequality, combined with a Minkowski type equality: 
$$|\partial\Omega|=\int_{\partial\Omega}p(s)H(s)\rd S,$$
where $p(s)=s\cdot \nu(s)$ is the support function of a star-shaped domain, $\nu(s)$ being the exterior normal. 

Using the normalized curvature flow, this Theorem has been improved in dimension $n=2$ to the more restrictive class of simply connected domain (note that spherical shells are simply connected only when $n\geq3$). One may ask wether the result holds true in any dimension among domains with connected boundary, but a counter example has been shown in \cite{FeNiTro16}: the authors construct a family of nodoids, diffeomorphic to a ball, in dimension 3, whose mean curvature is arbitrarily small, with a fixed volume. 

As a consequence of Theorem \ref{T:optimMean}, we get: 
\begin{corollary}
Let $\Omega$ be a domain which is simply connected if $n=2$, or star shaped if $n \geq 3$. Assume that $\Omega$ is not a ball and let $B$ be a ball of same volume. Then, for all $j\in \N$, there exists $\alpha_{0}>0$ such that 
$$\forall \alpha \in (\alpha_{0},+\infty), \quad \lambda_{j}^{R}(\alpha,\Omega)<\lambda_{j}^{R}(\alpha,B).$$
\end{corollary}

\chapter{The semi-classical magnetic Laplacian in corner domains}
\label{C:magnetic}
In this chapter we consider the semi-classical magnetic Laplacian 
\bel{E:magnlaplacian}
(-ih\nabla+A)^2=\sum_{j=1}^{n}(-ih\partial_{x_{j}}+A_{j})^2,
\ee
in a bounded simply connected domain $\Omega$, with magnetic Neumann boundary condition. All the results about the asymptotics of its first eigenvalue as $h\to0$ have the structure \eqref{E:megameta}. We show this asymptotics, with a remainder, when $\Omega$ is a three-dimensional corner domains. Unlike to the Robin laplacian (Section \ref{SS:Robincoin}), we are not able to perform a recursive procedure, and we have to proceed to an exhaustion of different model problems on the tangent geometries. We give various improvement of the remainder under stronger geometric assumptions, and we obtain a full asymptotic expansion when the domain is a lens, under some hypotheses. 

\section{The semi-classical magnetic Laplacian}
In dimension $n\in \{2,3\}$, the magnetic Laplacian with magnetic field takes the form \eqref{E:magnlaplacian}
where the vector field $A=(A_{1},\ldots,A_{n})$ is the magnetic potential. The associated magnetic field is the 2-form $B=\rd w_{\alpha}$, where $w_{A}:=\sum_{j=1}^{n}A_{j}\rd x_{j}$. In low dimensions, the magnetic field can be identified with $B=\curl A$. The operator is defined by the above differential expression, with a magnetic Neumann boundary condition $(-ih\nabla u-A)u \cdot n=0$. It is the Friedrichs extension of the form 
\[ H^{1}(\Omega)\ni u \mapsto  \int_{\Omega}|(-ih\nabla -A)u|^2,\]
\Bk
and an eigenpair $(\lambda,\psi)$ of this operator solves the boundary value problem
\begin{equation}
\label{eq:pbh}
   \begin{cases}
   (-ih\nabla+A)^2\psi=\lambda \psi\ \ &\mbox{in}\ \ \Omega \,,\\
   (-ih\nabla+A)\psi\cdot \nu=0\ \ &\mbox{on}\ \ \partial\Omega\,.
   \end{cases}
\end{equation}
\Bk

For a simply connected domain $\Omega$, due to the gauge invariance, the spectrum depends only on the magnetic field $B=\curl A$, and not on the choice of the magnetic potential. In this chapter, we denote by $\lambda_{j}(h)$ the $j$-th eigenvalue of the magnetic Laplacian.

The asymptotics of $\lambda_{j}(h)$ has receive a lot of attention for more than twenty years. One of the main motivations comes from the modeling of surface superconductivity, which involves the magnetic Laplacian with a large magnetic field, and can be linked to the semi-classical Laplacian described above (see  \cite{FouHel10} and the references therein).

From now on, we consider that $B$ is fixed. We assume that it is smooth enough  and, unless otherwise mentioned, does not vanish on $\overline\Omega$. 
The question of the semiclassical behavior of $\lambda_{1}(h)$ has been considered in many papers for a variety of domains, with constant or variable magnetic fields: Smooth domains \cite{BauPhiTa98, BeSt98,LuPan99-2, HeMo01, DpFeSt00, FouHe06, Ray09} and polygons \cite{Ja01,Bon06, BonDau06, BoDauMaVial07} in dimension $n=2$, and smooth domains \cite{LuPan00, HeMo02, HeMo04, Ray3d09, FouHel10} in dimension $n=3$. We will give an overview in the next paragraph, and we refer to \cite{FouHel10} and \cite{RayBook} for a more complete state of the art. Three-dimensional non-smooth domains were only addressed in two particular configurations---rectangular cuboids \cite{Pan02} and lenses \cite[Chap. 8]{Popoff} and \cite{PofRay13}, with special orientations of the magnetic field (that is supposed to be constant). Finally, in \cite{BoDauPof14}, general three-dimensional corner domains are treated, with some extensions to higher dimensions.

 Before describing this literature, we replace the problem in our framework by defining the tangent operators and their ground state energy.
 
 \begin{notation}
 \label{N:magnetic}
In this chapter, $L_{h}$ is the magnetic Laplacian associated with \eqref{eq:pbh}. For a cone $\Pi$ and a constant vector field $B_{0}\neq0$, we define $L_{h}[\Pi,B_{0}]$, the Neumann-magnetic Laplacian on $\Pi$ with a linear potential generating $B_{0}$, and $L[\Pi,B_{0}]:=L_{1}[\Pi,B_{0}]$. We denote by $E(\Pi,B_{0})$ the bottom of the spectrum of $L[\Pi,B_{0}]$.
For a singular chain $\dx=(x_{0},\ldots)\in \gC(\Omega)$, the local energy is defined by $E(\dx)=E(\Pi_{\dx},B(x_{0}))$, and, as in \eqref{D:lle}, $\sE=\sE(\Omega,B)$ is the infimum of $E(\cdot)$. Notice that $E(\dx)$ depends not only on the geometry of $\Pi_{\dx}$, but also on the value of the magnetic field at $x_{0}$. 
\end{notation}
Before \cite{BoDauPof14}, all the results about $\lambda_{1}(h)$ have the structure of \eqref{E:megameta}
with various estimates on $\lambda_{1}(h)-h\sE(B,\Omega)$, depending on the geometry, mostly on the form 
\bel{E:metaestimate}
-Ch^{\kappa_{-}} \leq \lambda_{1}(h)-h\sE(B,\Omega) \leq Ch^{\kappa_{+}}.
\ee
where $\kappa_{\pm}\in (1,2)$.

We describe now what is known on the values of the local energy, linking it with what is known for semi-classical asymptotics. We refer to \cite{FouHel10}, and more recently to \cite[Section 2]{BoDauPof14}, for a more detailed state of the art. 
 
 Using another scaling, one easily see that $E(\Pi,B)=|B|E(\Pi,\frac{B}{|B|})$, therefore we consider below the local energies for constant unitary magnetic field only.

The case where $\Pi=\R^{n}$ is explicit: here $E(\Pi,B)=1$ corresponds to the first Landau level.

 In dimension $n=2$, $B$ is just a scalar field, and $E(\Pi,1)=:\Theta_{0}\approx0.59$ when $\Pi$ is a half-plane. This value is the infimum of the eigenvalues of a one-dimensional family of Sturm-Liouville operators which will also appear in Part \ref{P:fibered}. For a sector $S_{\alpha}$ of opening $\alpha\in (0,2\pi)$, there holds $E(S_{\alpha},1) \leq \Theta_{0}$, with a strict inequality if $\alpha\in (0,\alpha_{0})$, were $\alpha_{0}$ is slightly above $\frac{\pi}{2}$, see \cite{Ja01,TBon,Bon06,ExLoPe17}. It is conjecture that the strict inequality holds if and only if $\alpha<\pi$, see \cite{BoDauMaVial07} for finite element computations by Galerkin projection.

Therefore, if $\Omega\subset\R^{2}$, and $b$ is a regular non-vanishing magnetic field
$$\sE(\Omega,B)=\min(\inf_{\Omega} B,\Theta_{0}\inf_{\partial \Omega}B,\min_{v\in \gA_{0}}(B(v)E(S_{\alpha(v)},1)),$$
where $\alpha(v)$ is the opening angle of a vertex $v\in \gA_{0}$. In that case, \eqref{E:metaestimate} is true with various values of $\kappa_{\pm}$, depending on the situation: 
 for the regular case, we refer to \cite{LuPan99, BeSt98,HeMo01,FouHe06} (constant magnetic field) and \cite{HeMo01,Ray09} (variable magnetic fields). For polygonal domains,  this is done in \cite{Bon06,BonDau06}.

 For regular domains in dimension 3, the local energy depends on the geometry as follows: denote by $\theta\in [0,\frac{\pi}{2}]$ the non-oriented angle between the  magnetic field and the boundary of the half-space $\Pi$, then, $B$ being unitary, $E(\Pi,B)$ depends only on $\theta$, and is denoted by $\sigma(\theta)$. Then, as proved in \cite{LuPan00,HeMo04}, $[0,\frac{\pi}{2}]\ni\theta\mapsto \sigma(\theta)$ is $C^{1}$ and increasing, moreover $\sigma(\frac{\pi}{2})=1$, and $\sigma(0)=\Theta_{0}$. This function is studied in more details as $\theta\to0$ in \cite{BoDauPopRay12}.
 Therefore, if $\Omega\subset\R^{3}$ is regular, 
$$\sE(\Omega,B)=\min(\inf_{\Omega} |B|,\inf_{\partial \Omega}\sigma(\theta(x))|B(x)|),$$
where $\theta(x)$ denotes the angle between $B$ and $\partial\Omega$ at  a point $x$. In the case where $B$ is constant, the minimum is reached at points on the boundary at which the magnetic field is tangent. Then, in \cite{HeMo04}, \eqref{E:metaestimate} is proved with $\kappa_{\pm}=\frac{4}{3}$. Under more geometrical hypotheses, a two term asymptotic is provided. The proof of \eqref{E:metaestimate} for a variable magnetic field is done in \cite{FouHel10} (with a remainder), and a more precise asymptotics is given in \cite{Ray3d09} under more hypotheses.

Concerning singular domains, rectangular cuboids has been addressed in \cite{Pan02}: \eqref{E:megameta} holds, and improvements are proved for particular configuration of the magnetic field. The case where $\Pi$ is a three-dimensional wedge (the model problem associated with an edge), is described in \cite{Pof13,Pof15}. As an application the situation where the domain $\Omega$ is a lens (two regular faces separated by a loop contained in a plane) is treated with a constant magnetic field orthogonal to the plane of the loop: \eqref{E:metaestimate} is obtained in \cite{Popoff} for lenses of small opening, with $\kappa_{\pm}=\frac{5}{4}$, and precised with a sharp asymptotics in \cite{PofRay13} under some non-degeneracy hypotheses, with $\kappa_{\pm}=\frac{3}{2}$.

\section{The need of a taxonomy in corner domains}
\label{S:taxo}
In the works described above, $\sE(\Omega,B)$ is known {\it de facto} or by hypothesis: for example, for a unitary magnetic field, it is $\Theta_{0}$ in the case of a regular domain or $E(S_{\frac{\pi}{2}},1)$ in the case of a square. But a major problem occurs in the case of a general corner domain: the minimizing tangent geometry is hard to determine. It is not even known for 2d corners, in which it is often assumed that the minimizing energy comes from the corner of smallest opening, a natural, although hard to check, hypothesis. 
 Our first step is to use Theorem \ref{T:Fsci} to prove the existence of a minimizer of the local energy. As for the Robin case, this step is important because it will show that $\sE(\Omega,B)>0$, and also because it will help us to construct a test function near a minimizer.

The monotonicity of $\dx \mapsto E(\dx)$ cannot be proved recursively as for the Robin Laplacian, because no analogous of \eqref{E:decomtensor} holds for the magnetic Laplacian. This property will come from the exhaustion of models occurring. In this paragraph, the magnetic field is constant and unitary. We use the second energy level $E^{*}$, defined similarly to \eqref{D:Estar0}. We have proved in \cite{BoDauPof14}:
\begin{theorem}
\label{T:dicho}
Let $\Pi\in \gP_{3}$, and $B\neq0$ a constant magnetic field. Then 
\bel{E:Monotochain}
E(\Pi,B) \leq E^{*}(\Pi,B).
\ee
Moreover, we have the following alternative: 
\begin{enumerate}
\item If $E(\Pi,B) < E^{*}(\Pi,B)$, then there exist an $L^{\infty}$ eigenvector for $L[\Pi,B]$ associated with $E(\Pi,B)$. 
\item If $E(\Pi,B) = E^{*}(\Pi,B)$, then there exists a singular chain $\dx \in \gC(\Pi)$ such that $E(\Pi_{\dx},B) < E^{*}(\Pi_{\dx},B)$ and $E(\Pi,B)=E(\Pi_{\dx},B)$.
\end{enumerate}
\end{theorem}
This theorem, stated in \cite[Section 7]{BoDauPof14}, relies on a taxonomy of model problems, since it is proved by the exhaustion of the different cones, according to their degree of singularity $d$:

For $d =1$, the cone $\Pi$ is a half-space, $E(\Pi,B)=\sigma(\theta)$ where $\theta$ is the angle between $B$ and $\partial\Pi$, and $E^{*}(\Pi,B)=E(\R^{3},B)=1$. Monotonicity of the local energy comes from $\sigma(\theta)\leq1$ and regularity of $\sigma$ was already known. 

Assume that $d=2$, i.e. that $\Pi$ is a wedge. Then $E^{*}(\Pi,B)$ can be expressed with $\sigma(\cdot)$, and the angles between the magnetic field and the two faces. The theorem comes from \cite{Pof15}. It relies on the fiber decomposition by partial Fourier transform along the edge, and the study of the associated band functions. In particular it is shown that the limit of these band functions are linked to $E^{*}(\Pi,B)$. Sufficient geometric conditions for the strict inequality are provided : this is true when the opening angle is small. Cases of equality are also provided in \cite{Pof13}.

 When $d=3$, it is proved in \cite{BoDauPof14} that the bottom of the essential spectrum of $L[\Pi,B]$ is $E^{*}(\Pi,B)$, leading to the inequality. In case of equality, the existence of a singular chain satisfying the theorem comes from application of Theorem \ref{T:Fsci} to the section of $\Pi$, together with the exhaustion of cases $d\leq2$. 
 Note that in case of strict inequality, the associated eigenfunctions have exponential decay. Examples of such cones can be found in \cite{Pan02,BoRay14,BoDaPoRa15}. 

As a consequence of Theorem \ref{T:Fsci} and \eqref{E:Monotochain}, we have: 
\begin{theorem}
Let $\Omega\in \gD(\R^{3})$, and $B\in C^{0}(\overline{\Omega})$ be a non-vanishing magnetic field. Then the function 
$$x\mapsto E(\Pi_{x},B(x))$$
is lower semi-continuous on $\overline{\Omega}$. Therefore, it reaches its infimum $\sE(\Omega,B)$, and we have that $\sE(\Omega,B)>0$.
\end{theorem}

\section{Asymptotics with remainder}

In this section, we describe our results about the asymptotics of $\lambda_{1}(h)$, as $h\to0$. Our main asymptotic result, proved in \cite{BoDauPof14}, is
\begin{theorem}
Let $\Omega\in \gD(\R^{3})$, and $A \in C^{2}(\overline{\Omega})$. Then there exists $C_{\Omega}>0$ and $h_{0}>0$ such that, for all $h\in (0,h_{0})$:
\begin{equation}
\label{eq:conv}
   \big|\ \lambda_{1}(h) - h \sE(\Omega,B) \big|\le
   \begin{cases}
   C_\Omega \big(1+ \|A\|^2_{W^{2,\infty}(\Omega)}\big) \, h^{11/10} ,
    &\mbox{$\Omega$ corner domain,} 
   \\[0.5ex]
   C_\Omega \big(1+ \|A\|^2_{W^{2,\infty}(\Omega)}\big) \, h^{5/4} ,
    &\mbox{$\Omega$ polyhedral domain.}
   \end{cases}
\end{equation}
Here the constant $C_\Omega$ only depends on the domain $\Omega$ (and not on $A$, nor on $h$).
\end{theorem}

Globally, the strategy of proof is the same that for the Robin Laplacian (and has been developed chronologically before): the lower bound comes from a suitable partition of the domain and approximation of the operator by local models, whereas the upper bound comes from the construction of test functions localized near a minimizer $x_{0}$ of the local energy. The test functions are provided by Theorem \ref{T:dicho}, and are then adapted to the local geometry, in the same way as for the procedure for the Robin Laplacian: in a row, they are qualified as sitting in case 1 of Theorem \ref{T:dicho}, and are centered at $x_{0}$, and as sliding in case 2, where they are issued from a higher singular chain than $(x_{0})$, their support is close to $x_{0}$ but avoids this point.

We describe now the differences, both in the proof and in the result, with the low-lying asymptotics for the Robin Laplacian.

Our result is valid in dimension 3, and the method does not extend in higher dimension. One of the main reasons is that the analysis of the local energy relies on a taxonomy of model problems. The analysis on wedges (\cite{Pof15}) is done by a fiber decomposition through a partial Fourier transform along the edge, and leads to a family of magnetic Laplacians, with a potential depending of the parameter. When the first band function reaches its infimum, it provides localization in the frequency dual to the edge, but the existence of such a minimum is not always guaranteed. In that case, it is still possible to link the limit of the band function to the second energy level $E^{*}(\Pi,B)$.
The specificity of this procedure does not seem to work directly in order to compare the bottom of the spectrum of a tangent operator on a cone $E(\Pi,B)$ and the second energy level $E^{*}(\Pi,B)$ for a general cone $\Pi$. This shows that the method proving the lower semi-continuity of the local energy does not adapt directly for the magnetic Laplacian, and it is not clear whether $\sE(\Omega,B)$ is not zero in higher dimensions. Nevertheless, we show in \cite{BoDauPof14} that for $\Omega \in \gD(\R^{n})$
$$\lim_{h\to}\frac{\lambda_{1}(h)}{h} = \sE(B,\Omega),$$
moreover we are able to prove a lower bound in any dimension, see \cite[Section 5.3]{BoDauPof14}, in the sense that for $\Omega\in \gD(\R^{n})$, there exist $\kappa\in (1,2)$ and $C>0$ such that for $h$ small:
$$\lambda_{1}(h) > h\sE(B,\Omega)-Ch^{\kappa},$$
moreover if $\Omega$ is polyhedral, $\kappa=\frac{5}{4}$. 

Another difference with the Robin Laplacian is the influence of the magnetic field, both in the minimization of the local energy and in the remainder estimates. In the analysis of Section \ref{SS:Robincoin}, the minimum of the local energy is given by the geometry of $\Omega$ only. In contrast here it may depends also on the variation of the magnetic field, giving rise to several different configurations for the locus of the concentration of the eigenfunctions. Therefore, it is not clear what could be an effective Hamiltonian for the next term in the asymptotics without additional hypotheses. 

 The error terms $R$ in \eqref{E:IMSschema}, which appears also in the upper bound, involve linearization of the metrics in $\cO(\eps(p)\kappa(p))=h^{\delta_{1}}$, and of the magnetic potential, in a combination of $\cO(\eta h )$ and $\cO(\eta^{-1}\eps(p)^4\kappa(p)^2)$, where $\eta>0$ is a parameter to be optimized. This is done by choosing $\eta=h^{\delta_{0}+2\delta_{1}-\frac{1}{2}}$. At the end, the optimization of scales leads to $\delta_{0}=\frac{3}{10}$ and $\delta_{1}=\frac{3}{20}$. In the case of a polygonal domain, there holds $\kappa=O(1)$, a one-scale analysis is sufficient and leads to $\delta_{0}=\frac{3}{8}$, as in \cite{HeMo01}.

In \cite{BoDauPof14}, we improve the upper bound of our asymptotics under stronger conditions: 

\begin{theorem}
\label{T:refinedasympt}
let $\Omega\in \gD(\R^{3})$ and $B\in C^{2}(\overline{\Omega})$ a magnetic field. 
\begin{enumerate}
\item
If $B$ vanishes somewhere in $\overline\Omega$, the lowest local energy $\sE(\Omega,B)$ is zero, and 
\begin{equation*}
  \lambda_{1}(h) \le    C_\Omega \big(1+ \|A\|^2_{W^{2,\infty}(\Omega)}\big) \, h^{4/3}.
\end{equation*}
\item
If there exists a corner $x_{0}\in \gA_{3}(\Omega)$ such that $\sE(\Omega,B)=E(\Pi_{x_{0}},B_{x_{0}})<E^{*}(\Pi_{x_{0}},B_{x_{0}})$, then 
$$
\lambda_{1}(h) \leq h \sE(\Omega,B) 
   +C_\Omega \big(1+ \|A\|^2_{W^{2,\infty}(\Omega)}\big) \, h^{3/2}|\log h| ,
   $$
   \item Assume that $\Omega$ is a straight polyhedron and $B$ is constant. Then 
   $$
\lambda_{1}(h) \leq h \sE(\Omega,B) 
   +C_\Omega \, h^{2} ,
   $$
\item
Assume that $A \in C^{3}(\overline{\Omega})$, then
\begin{equation}
\label{eq:convBIS}
    \lambda_{1}(h) - h \sE(\Omega,B) \le
   \begin{cases}
   C_\Omega \big(1+ \|A\|^2_{W^{3,\infty}(\Omega)}\big) \, h^{9/8} ,
    &\mbox{$\Omega$ corner domain,} 
   \\[0.5ex]
   C_\Omega \big(1+ \|A\|^2_{W^{3,\infty}(\Omega)}\big) \, h^{4/3} ,
    &\mbox{$\Omega$ polyhedral domain.}
   \end{cases}
\end{equation}
\end{enumerate}

\end{theorem}

The three first improvements are rather expected, and do not require new tool in the analysis. They come from a better optimization of error terms, using the hypothesis. But the asymptotics \eqref{eq:convBIS} requires a much deeper analysis of all the model problems described in Section \ref{S:taxo}, in order to make some terms vanish when evaluating the test functions.  The core of this upper bound is to take into account the fact that the functions coming from Theorem \ref{T:dicho}, used in the construction of our test functions, are not only $L^{\infty}$, but may have exponential decay in some direction. This is combined with subtle properties, such as that the Holder regularity of the local energy for wedges proved in \cite{Pof15}, and a good choice of gauge. 

Computing the next term in the asymptotics requires more hypotheses, but it not clear how stronger hypotheses on the geometry will turn into good properties for the minimum of the local energy. It is natural to look for a well (a punctual minimum) of the local energy. 

In a two-dimensional polygon, if the local energy is minimum at a corner (a natural hypothesis), with a gap (case 2 of Theorem \ref{T:refinedasympt}), then it is possible to give a full asymptotic expansion of the first eigenvalues: as described for the Robin Laplacian, if the model problem at each corner $v_{j}$ possesses $K_{j}$ eigenvalues below the bottom of the essential spectrum, it is possible to give the asymptotics of the first $K=\sum_{K_{j}}$ eigenvalues in power of $h^{1/2}$. We refer to \cite{Bon06,BonDau06}. These cases are called {\it corner concentration}.

For a regular domain in dimension 2, the local energy is nothing else but $B(x)$ in the interior and $\Theta_{0}B(x)$ on the boundary. If the minimum of the local energy corresponds to one (and only one) of the minimum of these two functions, in the sense that 
$$\min_{\Omega}B < \Theta_{0} \min_{\partial\Omega}B \ \ \mbox{or} \ \  \min_{\Omega}B > \Theta_{0} \min_{\partial\Omega}B,$$
and if this minimum is non-degenerate, reached in a unique point $x_{0}$, then the asymptotics of all the low-lying eigenvalues is available, using the concentration of the eigenfunctions near the point $x_{0}$. 

The case where $x_{0}$ is in the interior (\cite{LuPan00}) is in fact very close to the case without boundary (\cite{HeMo96,HeKo11}), and one obtains a full asymptotic expansion of all the low-lying eigenvalues. This case easily extends to higher dimension. 

When $x_{0}$ is on the boundary, it is also possible to localize (in phase-space) the operator near $(x_{0},k_{0})$, where $k_{0}$ is the frequency at wich the first band function of the model problem on $\R_{+}^2$ reaches its infimum $\Theta_{0}$, and to compute the full asymptotics expansions of $\lambda_{n}(h)$, see \cite{Ray13}.
In dimension 3, the situation is slightly more complicated because the local energy on the boundary is now $|B(x)|\sigma(\theta(x))$, where $\theta(x)$ is the angle between $B$ and the boundary at a point $x$. In the case where the minimum of the local energy corresponds to a unique non-degenerate minimum of this function, a three-term asymptotics is computed in \cite{Ray3d09}. The method should extend to the case where the local energy has a non-degenerate global minimum inside a given stratum.
 
The case of a constant magnetic field is very different in spirit because the local energy is now constant on a submanifold of a stratum: the whole boundary in dimension 2, and a curve inside the boundary in dimension 3. In dimension 2, a full asymptotic of the eigenvalues is available under the assumption that the curvature admits a unique non-degenerate minimum (\cite{FouHe06}). The method combines the localization of the eigenfunctions in phase space and a Grushin method. In dimension 3, only the first two terms of the asymptotics of the first eigenvalue are known, see \cite{HeMo04}. But these two results requiere a non-degeneracy hypothesis on the curvature, and the existence of an effective Hamiltonian leading the asymptotics, in the spirit of \cite{PankPof16}, is still an open question. The construction of a WKB expansion near the boundary is made in \cite{BoHerRay16}, and a suggestion for an effective operator on the boundary can be found in \cite{BoHerRay17}, for a two dimensional regular domain with a constant magnetic field. 

In \cite{PofRay13}, we consider a lens: a domain in $\R^{3}$ with two faces separated by an edge, which is a loop contained in the plane $\{x_{3}=0\}$. The magnetic field is $B=(0,0,1)$. For this domain, the local energy on the edge depends only on the opening angle $\alpha$ between the two faces. Moreover, it is increasing with respect to $\alpha$. Therefore, if this opening angle has a minimum $\alpha_{0}$, it will provide a minimum for the local energy. We assume that the minimum of the opening angle is non-degenerate, which implies the same property for the local energy. At this point, the model problem is a magnetic Laplacian on a 3d-wedge. We assume that for this model problem, $E<E^{*}$ (case 1 of Theorem \ref{T:dicho}), a fact which is check numerically and proved for $\alpha$ smaller than $\alpha_{m}\approx 0.38\pi$ (\cite{Popoff}). We finally need to assume that the first band function of the model problem has a unique non-degenerate minimum. Then we show under these assumptions the full asymptotic expansion 
$$\lambda_{j}(h)=h\sum_{p\geq0}c_{p,j}h^{p/4},$$
where $c_{0,j}=\sE$, $c_{1,j}=0$ and $c_{2,j}=\kappa_{0}+\omega_{0}(2j-1)$, where $(\kappa_{0},\omega_{0})$ are two constants. We also have an expansion of the eigenfunctions. Once again, our tools are based on Agmon estimates for the localization in phase and space of the eigenfunctions, and a Grushin procedure in order to approach the magnetic Laplacian by an effective operator with a potential having a non-degenerate minimum at the low-lying energies. 

\part{Translation invariant magnetic fields (and their perturbations)}
\label{P:fibered}
\chapter[Invariant magnetic fields]{Spectrum of magnetic Laplacians with translation invariant magnetic field}
In this chapter, we present the {\it free operators} of this part: they are magnetic Laplacians with magnetic fields having translation invariant properties. More precisely, they have the form \eqref{E:translationnally}. In section \ref{S:threshold}, we define precisely the systems we will consider, and we give their fiber decomposition through partial Fourier transform. The associated fiber operators have discrete spectrum, the so called {\it band functions}. We present the notion of threshold of a fibered operator from \cite{GeNi98}, and explain that our operators have a new kind of threshold, corresponding to a finite limit of a band function. In Section \ref{S:bandglobal}, we give their asymptotics for large frequencies. In Section \ref{S:current}, we exploit these results by describing the {\it bulk states}: they are quantum particles localized in energy near a threshold, with very small propagation along the direction of invariance.

\section{The thresholds of fibered operators}
\label{S:threshold}
In this section, we defined in details the different Schr\"odinger operators we will analyze, and we discuss the notion of thresholds in their spectrum.

The operators considered in the two following parts have the form \eqref{E:translationnally}, (i.e. $(-i\nabla-A)^2$, acting on $L^{2}(\omega\times \R)$), more precisely, we will mainly considered the following classes of Hamiltonian: 
\begin{enumerate}[label=\Alph*]
\item\label{A} Half-plane submitted to a constant magnetic field. Here $\omega=\R^{_{+}}$ and $B=1$. The operator is defined by a boundary condition at $x=0$ (mainly, Neumann or Dirichlet). This model has been studied in numerous model arising in surface superconductivity and quantum Hall effect. We refer to for \cite{SaGe63,dBiePu99,MaMaPu99,HeMo01,BucGra06,HisSoc08,FouHel10,BruMirRai13}) for a non exhaustive list of articles where this operator is involved. \label{popoff}
\item\label{B} Iwatsuka magnetic field. Here, $\omega=\R$ and $B=B(x)$ is increasing and have finite limits. Moreover, two cases are particularly relevant: 
\begin{enumerate}[label=\arabic*]
\item\label{B1} Standard Iwatsuka magnetic field. This is the case where $B$ is positive regular with positive limits. This model has been studied first as an example of a magnetic Laplacian with absolutely continuous spectrum (\cite{Iwa85}), and present transport properties along the $y$ direction (\cite{PeetVas93,ManPur97,HisSoc14}).
\item\label{B2} Magnetic steps. This is the case where $B$ is constant, with values $B^{\pm}$, on each half-line $\R^{\pm}$. This model presents interesting transport properties in the both directions along the $y$ axis, and can be solved analytically (\cite{PeetRej00,HisPerPofRay16}), it is also involved in surface superconductivity with a discontinuous constant magnetic field (\cite{AssKach16}).
\end{enumerate}
\item\label{C} Magnetic {\it wire}. Here, $\omega=\R^{2}$, and $B=b(r)(\sin \theta,-\cos \theta,0)$, where $(r,\theta,z)$ are cylindrical coordinates of $\R^3$. These models have been studied in \cite{Yaf03,Rai08,Yaf08}. In that case, the magnetic potential can be chosen to be $A(r)=(0,0,a(r))$. The case $a(r)=\log r$ corresponds to a magnetic field created by an electric wire placed in the $z$ axis (\cite{Yaf03,BruPof15}), the case $a(r)=r$ to a unitary magnetic field (\cite{Popoff,Gen18}) and more general cases are treated in \cite{Yaf08}.
\end{enumerate}
These models can be fibered through the partial Fourier transform $\cF$ along the direction of invariance, and satisfy \eqref{E:fibering}. The first two models lead to the same one dimensional operator 
\bel{E:ghk}
\gh(k)=-\partial_{x}^2+(a(x)-k)^2, \quad x\in \omega
\ee
where $a'(x)=B(x)$, and with a boundary condition (Neumann or Dirichlet) when $\omega=\R^{+}$, model \ref{A}. With our hypotheses, the domain of these operators does not depend on $k$. The fiber operators associated with model \ref{C} is slightly different, and the relevant case $a(r)=\log r$ will be described in Section \ref{S:magneticwire}.
 In all cases, except when $B$ has zero as a limit in case \ref{B2} (this case is treated in \cite{HisPerPofRay16}), the fiber operator $\gh(k)$ has compact resolvent, and the resolvent $(\gh(k)+i)^{-1}$ is analytic with respect to $k\in \R$. We denote by $\lambda_{j}(k)$ its eigenvalues, which are non-degenerate, and therefore analytic, and we denote by $u_{j}(\cdot,k)$ associated normalized eigenfunctions. In particular, since we have \eqref{E:fibering}, the spectrum is absolutely continuous, provided the band functions are not constant, which will be the case in our models. 

%
%

A definition of thresholds for fibered operators is given in \cite{GeNi98}. In this article, G\'erard and Nier define the class of {\it analytically fibered operators} (see \cite[Section 2]{GeNi98} for the precise definitions). They define the set of {\it energy-momentum} 
$$\Sigma :=\{ (\lambda,k)\in \R\times M, \lambda\in \sigma(\gh(k))\},$$
where the momentum space, $M$, a real analytic manifold, is $\R$ in our case. Then, they define the projection $p_{\R}:\Sigma\to \R$ by $p_{\R}(\lambda,k)=\lambda$. Assuming that this projection is proper, they show the existence of a stratification of the two sets $\Sigma$ and $\R$ associated with $p_{\R}$, in the sense that the dimensions of the strata agree with the rank of $p_{\R}$ on each stratum. The result relies on the theory of analytic singularities. Then, the thresholds are defined as the strata of dimension 0 of $\R$. They therefore form a discrete subset of the spectrum. The authors show a Mourre estimate outside this set, based on the construction of a conjugate operator. 

This theory has its roots in the analysis of the spectrum of the Schr\"odinger operator with periodic potential, which can be fibered through Floquet-Bloch transform. The energy-momentum set $\Sigma$ is known as the Bloch variety. In this framework, a quite complete analysis of the thresholds can be done, see \cite{Ge90}. In particular, when a threshold corresponds to a critical point $E$ of $\lambda_{j}(\cdot)$, a Taylor expansion of $\lambda_{j}$ at this point is available. It can provide a good description of the branching point of the resolvent near $E$. This point is also a key argument in the existence of an effective Hamiltonian, in the study of the number of eigenvalues in the gap of the spectrum for perturbations of such an operator, see \cite{Rai92}, and also the next chapter.

Our setting is slightly different. Firstly, this general framework is not necessary because the variable $k$ lies in $\R$, and the operator $\gh(k)$ acts in dimension 1 (except for the model studied in \cite{Pof15} where the dimension is 2), and the eigenvalues are simple, which involve that the band functions  are analytic. In that case standard thresholds correspond to the critical points of the band functions $\lambda_{j}(\cdot)$. The point making our operators out of the setting of \cite{GeNi98} is that the band functions are not proper in general, and may possess finite limits as $|k|\to+\infty$, giving rise to a new kind of threshold.

We define the thresholds for the above models as either one of the critical point of the band functions, or one of their limits.

Our goal is not to define thresholds in a general setting, but more to present results and methods for the study of the spectrum of $H_{0}$ (and its perturbations) near these values. An analytic proper band function has, at a critical points $k_{0}$, a Taylor expansion, which allows to remplace formally the operator by the symbol provided by its Taylor expansion, for frequencies $k \approx k_{0}$. In our cases, we will need to replace this idea by an asymptotic expansion of the band function as $k\to\infty$, and also to use the expansion of the corresponding eigenfunctions. As we will see, the band functions may present different behavior, depending on the model. 

\section{Asymptotic behavior of the band functions}
\label{S:bandglobal}
We collect here known and elementary properties of the band functions. The objective is to illustrate that their behaviors at infinity can be very different, depending on the model.

\begin{theorem}[\cite{DauHe93,HeMo01} for the global behavior, \cite{Ivrii,Popoff} for the asymptotics] Consider the model \ref{A}. In case of a Dirichlet boundary condition (D), the band functions are decreasing. In case of a Neumann one (N), they admit a unique non-degenerate minimum. In both cases, 
$$\lim_{k\to-\infty} \lambda_{j}(k)=+\infty \quad \mbox{and} \quad \lim_{k\to+\infty}\lambda_{j}(k)=2j-1:=E_{j},$$
moreover, $$\lambda_{j}^{D/N}(k)-E_{j}\underset{k\to+\infty}{\sim}\pm c_{j}k^{2j-1}e^{-k^2},\quad c_{j}>0.$$
Moreover the asymptotics of $\lambda_{j}'$ at $+\infty$ can be deduced by differentiating the above formula. 
\end{theorem}

\begin{theorem}
[\cite{Iwa85,ManPur97}
\label{T:asymptIwa} for the global behavior, \cite{MirPof18} for the asymptotics] Consider the model B\ref{B1}. Assume that the Iwatsuka magnetic field is increasing, positive, $C^{1}$, with $\lim_{x\to \pm\infty}=B_{\pm}$. Then the band functions are increasing, with 
$$\lim_{k\to \pm \infty}\lambda_{j}(k)=E_{j}B_{\pm}:=\cE_{j}^{\pm}.$$
Assume moreover that the magnetic field has a power-like convergence: 
\bel{E:PowerlikeB}
\exists x_{0}>0, \exists M>0, \forall x \geq x_{0}, \quad  B(x)=B_{+}-x^{-M}.
\ee
Then $$\lambda_{j}(k)-\cE_{j}^{+}\underset{k\to+\infty}{=}-\frac{E_{j}B_{+}^{M}}{k^{M}}+\cO(k^{-M-2}).$$
\end{theorem}
There exist a slightly more general version of this theorem, for a magnetic field which satisfies \eqref{E:PowerlikeB} only asymptotically, see \cite[Theorem 2.2]{MirPof18}. This theorem includes a two-term asymptotic, but we have preferred to present the above version, easier to read.

The case of a magnetic step, model B\ref{B2}, offers a different behavior, which depends on the ratio of the two values of the magnetic field (\cite{HisPerPofRay16}).  Using a normalization, we assume that $B_{+}=1$, and that $B_{-}\in [-1,0)$. We denote by $b=|B_{-}|$. Notice moreover that by a symmetry argument, the band functions of the case $B_{+}=-B_{-}=1$ is the collection of the band functions of the model in half-plane with both Neumann and Dirichlet boundary condition.

\begin{definition}\label{def}
We define $\mathfrak{L}=\left\{E_{j},j\geq 1\}\bigcup\{bE_{j}, j\geq 1\right\}$, and we define the splitting set $\mathfrak{S}=\left\{E_{j},j\geq 1\}\bigcap\{bE_{j}, j\geq 1\right\}$. We will say that there is a \emph{splitting of levels} if $\mathfrak{S}\neq\emptyset$ i.e. if there exist positive integers $\ell<j$ such that
\beq\label{eq:split1}
b=\frac{2\ell-1}{2j-1}.
\eeq
\end{definition}

The set of limit points of the band functions $\lambda_{j}(k)$ as $k \to +\infty$ is precisely $\mathfrak{L}$.
Moreover, we may describe the asymptotics of the band functions in the following way (\cite{HisPerPofRay16}): 
\begin{theorem}Consider the model \ref{B2}.
\begin{enumerate}
\item Non-splitting case.  Let $\lambda\in \mathfrak{L}\setminus\mathfrak{S}$. We can write in a unique way $\lambda=E_{j}$ or $\lambda=bE_{j}$ with $j\geq1$. 
Then there exists a unique $p_{j}\in \N^{*}$ such that:
\begin{enumerate}
\item $\lambda_{p_{j}}(k)-E_{j}=\eps_{j}(k)\underset{k\to+\infty}{\to} 0$ (in the case $\lambda=E_{j}$),
\item $\lambda_{p_{j}}(k)-bE_{j}=\eps_{j}(k)\underset{k\to+\infty}{\to} 0$ (in the case $\lambda=bE_{j}$),
\end{enumerate}
with, in the first case,
\[
\eps_{j}(k)=-\frac{1}{\sqrt{\pi}}2^{j-2}(1+b)\frac{1}{(j-1)!}k^{2j-3}e^{-k^2}(1+o(1)),
\]
and in the second case,
\[
 \eps_{j}(k)=-\frac{1}{\sqrt{\pi}}2^{j-2}b^{3/2-j}(1+b)\frac{1}{(j-1)!}k^{2j-3}e^{-k^2/b}(1+o(1)).
\]

\item Splitting case.  Assume that the splitting set is not empty and consider $\lambda\in\mathfrak{S}$. Let us consider $j\geq 1$ and $\ell\geq 1$ such that $\lambda=E_{j}=bE_{\ell}$. There exists $p_{j}\in \N^{*}$ such that
\[
\left\{
\begin{aligned}
&\lambda_{p_{j}}(k)-E_{j}=\eps_{j}^{-}(k)
\\
&\lambda_{p_{j}+1}(k)-E_{j}=\eps_{j}^{+}(k)
\end{aligned}
\right.
\]
with, as $k$ tends to $+\infty$,
\begin{align*}
\eps_{j}^{-}(k)&=-\frac{2^{j-3/2}(1+b)}{(j-1)!\sqrt{\pi}}k^{2j-3}e^{-k^2}(1+o(1)),
\\
\eps_{j}^{+}(k)&=\frac{2^{\ell+\frac{3}{2}}b^{-\ell+\frac{3}{2}}}{(\ell-1)!(1+b)\sqrt{\pi}}k^{2\ell+1}e^{-k^2/b}(1+o(1)).
\end{align*}
\end{enumerate}
\end{theorem}
\begin{remark}
These asymptotics permit to state when a band function is not monotonous, but there is no result about the number of minima. 
\end{remark}

These asymptotics are quite elementary to reach. The limit being the Landau levels can be understood as follows: as $k$ gets large, the eigenfunction associated with a bounded eigenvalues will concentrate in a zone $|a(x)-k|=\cO(1)$. But, as $k$ gets large, the value of the magnetic field in this zone is close to a constant (or even equal to a constant in the case of a magnetic step). These arguments show that the spectrum $(\lambda_{j}(k))_{n\geq1}$ converge the Landau levels. The computation of the next term needs a little more work. After some elementary changes of variable, the operator takes the form $h^{-1}\left( -h^2\partial_{X}+W(X,h)\right)$, where $W\geq0$ is confining and vanishes. For the Iwatsuka model \ref{B1}, $W$ has a unique minimum with an explicit Taylor expansion. The method, standard, consist in constructing quasi-modes, in the spirit of the harmonic approximation. This method tends to fail for magnetic field which are constant at infinity, as in the cases \ref{A} and \ref{B2}, because the potential $W$ will be equal to $X^2$ near its minima. In the case of a magnetic step, the potential $W$ is piecewise quadratic, and reaches its minimum in two points. The exponential remainder can then be understood as a tunneling between wells. But in those cases, the use of parabolic cylindrical functions solving the eigenvalue problem allows to recover the asymptotics.

\section{Study of currents}
\label{S:current}
%

In this section we present notions related with the quantum motion of a particle submitted to a translation invariant magnetic field. Firstly, we introduce the current in the direction of invariance, and express it using the group velocity (the derivative of the band function). The notion of edge states has appear in numerous models, such as in \cite{Ha82,ManPur97,ExJoyKov99,dBiePu99,PeetRej00,HisSoc14}, and roughly, corresponds to particle localized in energy far from the thresholds. Their counterpart, the so-called bulk states, are localized in energy near the thresholds, and we provide a method to have quantitative estimates of their current.

Recall that $\Omega=\omega \times \R$ is assumed to be translation invariant in the last variable, and denote by $\cF$ the partial Fourier transform in this variable. Given a function $f\in L^{2}(\Omega)$, we define its $j$-Fourier coefficient by 
$$f_{j}(k)=\langle (\cF f)(\cdot,k),u_{j}(\cdot,k) \rangle$$
and its $j$-th harmonic by the projection 
$$\pi_{j}(f)(x,y):=\frac{1}{\sqrt{2\pi}}\int_{\R} e^{iky}u_{j}(x,k)f_{j}(k) \rd k. $$


 The current operator along $y$ is 
 $$J_{y}=-i[H_{0},y],$$
and direct computations show that 
$J_{y}=-i\partial_{y}-a$. If we denote by $y(t):=e^{-itH_{0}}ye^{itH_{0}}$, the evolution of the position in the $y$ direction, then there holds $\frac{\rd y}{\rd t}(t)=e^{-itH_{0}}J_{y}e^{itH_{0}}$, so the velocity in the $y$ direction is the evolution of the observable $J_{y}$. 

This operator is naturally linked to the the derivative of the band functions, the so-called ``group velocity'', by the Feynman-Hellmann formula: 
$$\forall f \in \dom(H_{0}), \quad \langle J_{y}\pi_{j}f, \pi_{j}f\rangle=\int_{k\in \R}\lambda_{j}'(k)|f_{j}(k)|^2 \rd k.$$
Extending this formula to compute $\langle J_{y}f,f\rangle$ for any $f \in \dom(H_{0})$ is not direct, because interaction terms appear. However, if $f$ satisfies the {\it non-overlapping condition} $i\neq j \implies \supp(f_{i})\cap \supp(f_{j})=\emptyset$, then there holds 
$$\quad \langle J_{y}f,f\rangle=\sum_{j\geq1 } \int_{k\in \R}\lambda_{j}'(k)|f_{j}(k)|^2 \rd k.$$
Therefore, estimates on the current will be provided by the control of the derivative of the band functions. 

Assume that $\lambda_{j}$ is strictly monotone (if not, the following procedure can be easily adapted). We may assume that $\lambda_{j}$ is increasing, the case of a decreasing band function being the same. Let $\eta_{j}$ be the reciprocal function of $\lambda_{j}$. For an interval $I\subset \sigma(H_{0})$, if we denote by $\dP_{I}$ the spectral projection associated with $H_{0}$ on $I$, we say say that $f$ is localized in energy in $I$ when $\dP_{I}f=f$, i.e. when $f \in \ran(\dP_{I})$. In that case for all $j\geq1$, we have that $\supp(f_{j})\subset \eta_{j}(I)$, and we can control the current of each non zero harmonic $\pi_{j}f$ as follow: 
\bel{E:metacurrent}
\inf_{\eta_{j}(I)}\lambda_{j}' \leq \frac{\langle J_{y}\pi_{j}f, \pi_{j}f\rangle}{\| \pi_{j}f\|^2} \leq \sup_{\eta_{j}(I)}\lambda_{j}' 
\ee
We see now the role of the set of thresholds $\cT$, defined by the union of the critical points and of the limits of the band functions: if $\overline{I}\cap \cT=\emptyset$ and if a finite number of band function cross in $\overline{I}$, then there exists a constant $C>0$ (which does not depend on $f$) such that
$$ \forall j\geq1, \quad C \| \pi_{j}f \|^2 \leq \langle J_{y}\pi_{j}f, \pi_{j}f\rangle, $$
see \cite{dBiePu99} for the analysis in model \ref{A} and  \cite{ManPur97} for model \ref{B}.
The hypothesis that a finite number of band functions cross $I$ is not satisfied in magnetic-wire models \ref{C} (see \cite{BruPof16,Gen18}), where it is harder to get estimates on the current. We refer to \cite[Section 3]{Gen18} for an estimate of the current for a unitary magnetic field. 
 
If the energy interval $I$ contains a threshold, it is clear that no lower bound on the current is available (\cite{dBiePu99,HisSoc14}). The approach, developed in \cite{HisPofSoc14}, is to consider an energy interval close to a threshold and to use the asymptotic expansion of the band functions converging to this threshold. To make it simple let us assume that the band function $\lambda_{j}(k)$ converges to a threshold $\cE_{j}$ as $k\to+\infty$, and consider $I:=(\cE_{j}\pm\delta,\cE_{j})$ (the term $\pm \delta$ depends on wether the convergence is by below or by above), with $0<\delta \ll1$. Here, $\cE_{j}=E_{j}=2j-1$ in model \ref{A} and $\cE_{j}=B_{+}E_{j}$ in model \ref{B}. Then, using the asymptotics of Section \ref{S:bandglobal}, we get in \cite{HisPofSoc14} and \cite{MirPof18}:

\begin{theorem}
Consider the model \ref{A}, with Dirichlet boundary conditions. Let $j\geq1$ be fixed. There exist $\delta_{0}>0$ and $C>0$ such that $\forall \delta\in (0,\delta_{0})$,
$$\forall f\in \ran(\dP_{I}), \quad |\langle J_{y} \pi_{j}f,\pi_{j}f\rangle | \leq \left(2\delta \sqrt{|\log \delta|} + C\frac{\delta \log|\log \delta|}{\sqrt{|\log \delta|}}\right) \|\pi_{j}f\|^2.$$
Moreover, there exists another constant $c>0$, such that for all $\eps\in (0,1)$, there exists $\delta_{0}>0$ such that $\forall \delta\in (0,\delta_{0})$, 
$$\forall f\in \ran(\dP_{I}), \quad \int_{x=0}^{(1-\eps)\sqrt{|\log \delta|}}\int_{y\in\R} |\pi_{j}f(x,y)|^2\rd y \rd x \leq C \eps^{2j-1}\delta^{\eps^2}|\log \delta|^{\frac{2j-1}{2}}(1-\eps^2)\|\pi_{j}f\|^2 .$$
\end{theorem}

\begin{theorem}
Consider the model B\ref{B1}, with a magnetic field as in Theorem \ref{T:asymptIwa}. There exist $\delta_{0}>0$ and $C_{j}>0$ such that $\forall \delta\in (0,\delta_{0})$,
$$\forall f\in \ran(\dP_{I}), \quad |\langle J_{y} \pi_{j}f,\pi_{j}f\rangle | \leq \left(\beta\delta^{1+\frac{1}{M}}+ C\delta^{\frac{\eta}{M}}\right) \|\pi_{j}f\|^2$$
with $\beta=\frac{M}{B_{+}E_{j}^{\frac{1}{M}}}$ and $\eta=\min(2M+2,M+3)$.
Moreover, there exists another constant $c>0$ such that $\forall \delta\in (0,\delta_{0})$, 
$$\forall f\in \ran(\dP_{I}), \quad  \int_{-\infty}^{c\delta^{-\frac{1}{M}}}\int_{y\in\R}|\pi_{j}f(x,y)^2|\rd y \rd x \leq C \delta^{\frac{\eta}{M}}$$
\end{theorem}

In the above results, the estimates on the currents is a direct consequence of \eqref{E:metacurrent} and of the asymptotics on the derivative of the band functions. The result on the localization of functions localized in energy in $I$ needs more work. The condition $f\in \ran(\dP_{I})$ on the energy, together with $I$ close to a threshold, imply that $\supp(f_{j})$ is located in the zone $k \gg 1$. Here, a typical magnetic effect arise: using that the symbol $\xi_{x}^2+(k-a(x))^2$ of $H_{0}$ should lie in $I$ (when restricted to $\ran \dP_{I}$), we get that $\pi_{j}f$ should concentrated in the zone $a(x)-k=O(1)$, therefore in a zone corresponding to large values of $x$. These considerations become rigorous by using precise localization estimates for the eigenfunction $u_{j}(\cdot,k)$, as $k\to+\infty$.

\section{Decomposition of the free resolvent and absorption principle}
In this section we explain how to extend the resolvent near a point in the spectrum in a suitable functional space. The method is well known outside thresholds, using weighted spaces, and we describe functional spaces in which the resolvent is well defined up to the thresholds.

Let $R(z)=(H_{0}-z)^{-1}$ be defined on $\C \setminus \sigma(H_{0})$. For $f$ and $g$ in $L^{2}(\Omega)$, we have, using the partial Fourier transform $\cF$:
\bel{E:decomposR}
\langle R(z)f,g\rangle =\sum_{j \geq 1}\int_{k\in\R} \frac{f_{j}(k)\overline{g_{j}(k)}}{\lambda_{j}(k)-z}\rd k.
\ee
In this expression, we clearly see that one of the integrals may be undefined as $z=E\in \ran(\lambda_{j})$. One of the possible strategies to define the resolvent is to add conditions on $f$ and $g$ such the integral converges in a suitable sense. The case where $E$ is reached by a band function is well known, and the resolvent admits a limit as $z=E+it$ and $t\downarrow0$ (denoted by $z=E+i0^{+})$, in weighted spaces (\cite{Va66,CroDer95,GeNi98,So01}). We recall below the procedure by a direct construction, avoiding Mourre estimates. Assume that the energy $E$ correspond to a non-critical point of a band function $\lambda_{j}$, in the sense that 
$$\exists k_{0}\in \R, \quad \lambda_{j}(k_{0})=E, \lambda_{j}'(k_{0})\neq0,$$
then, assuming that $I$ is an open interval containing $k_{0}$ on which $\lambda_{j}$ is a diffeomorphism, we isolate the singular term in \eqref{E:decomposR}, an perform the change of variable $\lambda=\lambda_{j}(k)$:
\bel{E:intsing}
\int_{I}\frac{f_{j}(k)\overline{g_{j}(k)}}{\lambda_{j}(k)-z}\rd k=\int_{\lambda_{j}(I)}\frac{\tilde{f_{j}}(\lambda)\overline{\tilde{g_{j}}(\lambda)}}{\lambda-z} \eta_{j}'(\lambda) \rd \lambda,
\ee
where $\tilde{f}_{j}:=f_{j}\circ \eta_{j}$ (and similar notation for $g_{j}$), $\eta_{j}$ being the reciprocal function of $\lambda_{j}$ on $I$. This integral can be defined as $z= E+i0^{+}$ as a principal values when $\tilde{f_{j}}$ and $\tilde{g_{j}}$ are Lipschitz functions near $E$. Let $s>\frac{1}{2}$, a sufficient condition for that is that $f$ and $g$ are in the weighted space $L^{2,s}(\Omega):=\{f \in L^{2}(\Omega), \langle y \rangle^{s}f \in L^{2}(\Omega)\}$, see \cite{CroDer95}:

\begin{proposition}
Consider model \ref{A} or \ref{B}. Let $K$ be a compact subset of $\sigma(H_{0})$ without thresholds and $s>\frac{1}{2}$. Then $z\mapsto R(z)$, initially defined on $\C\setminus\sigma(H_{0})$, extends to $\overline{\C^{+}\cap K}$ as a H\"older function, for the topology of the norm in $\cB(L^{2,s}(\Omega),L^{2,-s}(\Omega))$.
\end{proposition}
The case where $k_{0}$ is a critical point of a band function can be treated by adding additional vanishing conditions for $f_{j}$ at $k_{0}$, see \cite{CroDer95,So01}. 

We now describe how to extend the resolvent when $E$ is the limit of one of the band functions. 
Assume that $E=\cE_{j}$ is the limit of exactly one band function $\lambda_{j}$, at $+\infty$, which is monotonous. The case where $\lambda_{j}$ is not monotonous (but monotonous in a neighborhood of $+\infty$) can be treated in the same way. To make it simple, assume finally that $E$ is not a threshold of another band function. Then, as above, for $z\in \C\setminus\sigma(H_{0}),$ we consider the $j$-th term in \eqref{E:decomposR}, and write it as in \eqref{E:intsing}. But this time, $E$ is an end point of the energy interval $\lambda_{j}(\R)$, moreover, $\eta_{j}'(\lambda)\to+\infty$ as $\lambda\to E_{j}$, therefore we need additional conditions to define the integral as a principal value.

Denote by $I_{j}:=\lambda_{j}(\R)$ (remark that $\cE_{j}$ is an endpoint of $I_{j}$), and $C^{0,\alpha}(I_{j})$ the space of the H\"older functions of $I_{j}$ endowed with its norm $\|\cdot\|_{C^{0,\alpha}(I_{n})}$. In what follows, the pair $(s,\alpha)$ satisfies
\begin{equation}
\label{E:salpha}
s>\frac{1}{2} \ \ \mbox{and} \ \ 0\leq\alpha<\min(1,s-\tfrac{1}{2}).
\end{equation}
Introduce the vector space
$$\cX^{s,\alpha}_{j}:=\{ f\in L^{2,s}(\Omega),\eta_{j}^{1/2}\tilde{f_{j}}\in C^{0,\alpha}(\overline{I_{j}})\cap L^{2}(I_{j}), \left(\eta_{j}^{1/2}\tilde{f_{j}}\right)(E_{j})=0\}, $$
endowed with its natural norm
$$\|f\|_{\cX^{s,\alpha}_{j}}:=\|\eta_{j}^{1/2}\tilde{f_{j}} \|_{C^{0,\alpha}(I_{j})}+\| \eta_{j}^{1/2}\tilde{f_{j}} \|_{L^2(I_{j})}+\|f\|_{L^{2,s}(\Omega)}.$$
We have, see \cite{PofSoc15}:
\begin{lemma}
The space $\cX^{s,\alpha}_{j}$ is a Banach space, dense in $L^{2}(\Omega)$.
\end{lemma}
For a compact set $K$, we define $\cT_{K}:=\{j\in \N, \cE_{j}\in K\}$, and 
$$\cX^{s,\alpha}_{K}:=\bigcap_{j\in \cT_{K}}\cX^{s,\alpha}_{j}.$$ Then: 
\begin{theorem}
\label{T:PALseuil}
Let $H_{0}$ be one of the Hamiltonians \ref{A} or B\ref{B1}, $K\subset \C$ a compact set, and $(s,\alpha)$ as in \eqref{E:salpha}. Then, $R(z)$, initially defined on $\C^{+}\cap K$, extends to an $\alpha$-H\"older continuous functions on $\overline{\C^{+}\cap K}$ in the uniform operator topology in $\cB(\cX^{s,\alpha}_{K},(\cX^{s,\alpha}_{K})')$
\end{theorem}
This result has been written for the model \ref{A} in \cite{PofSoc15}. Let us now describe some properties of this absorption space. One of the conditions provides
$$f\in \cX^{s,\alpha}_{j}(\Omega) \implies w_{j}^{\alpha}f_{j}\in L^{\infty}(\R) $$
with $w_{j}(k):=|\lambda_{j}(k)-\cE_{j}|^{-\alpha}|\lambda_{j}'(k)|^{-\frac{1}{2}}$. Therefore, functions in the absorption space have decay properties for its $j$-Fourier coefficient. This has several consequences, such as the regularity in the dual variable $y$, as expected. Another, more surprising property, is the decay in the $x$ variable. Once again, this is a magnetic effect, arising from the mixing in phase space of the variables $k$ and $x$ through the symbol of $H_{0}$. We illustrate these two properties on the model \ref{A} with Dirichlet boundary condition, for which $w_{j}(k)\underset{k\to+\infty}{=}C_{j,\alpha}k^{\alpha-j(2\alpha+1)}e^{k^2(\alpha+\frac{1}{2})}$, by the result below (see \cite[Section 3]{PofSoc15} for the detailed proof): 

\begin{theorem}
Consider model \ref{A}. Under the assumptions of Theorem \ref{T:PALseuil}, for any positive real number $\beta < \min \left( 1,\frac{2\alpha+1}{1+\sqrt{2\alpha+1}} \right)$, we may find two constants $C_j(\beta)>0$ and $L_{j}(\beta)>0$ such that we have:
\bel{g0}
\forall f \in \cX_{j}^{s,\alpha}(\Omega),\ \forall L \geq L_j(\beta),\quad \int_{L}^{+\infty} \| \pi_{j}f(x,\cdot) \|_{L^2(\R)}^2 \rd x \leq C_j(\beta) \| f \|_{\cX_{j}^{s,\alpha}(\Omega)}^2 e^{- \beta L^2}.
\ee 
Moreover, the function 
$y \mapsto \|\pi_{j}f(\cdot,y)\|_{L^2(\R_{+})}$, initially defined in $\R$, admits an analytic continuation to $\overline{\C^{-}}$.

\end{theorem}

\chapter{Effect of perturbations on thresholds}

In this chapter we will consider a fibered Hamiltonian $H_{0}$, such as those of the previous chapter, and perturb it by a sign-definite potential $V$. For suitable perturbations, the essential spectrum will remain the same, but discrete eigenvalues may appear. It is well known that near the thresholds of the spectrum of $H_{0}$, an accumulation of eigenvalues can occur. This phenomenon is deeply related to the nature of the threshold. 

%

The case of a constant (homogeneous) magnetic field has been widely studied, starting from \cite{AvrHerSi78,Sobolev}, where the most striking behaviour is the fact that even a compactly supported perturbation can create an infinite number of eigenvalues. This framework, where the band functions (and the associated eigenfunctions) are more explicit, has then received a lot of attention, see \cite{Tam88,Rai90,Ivrii0,RaiWar02,BonBruRai07}, and has been extended to various kind of perturbations.
 
In \cite{Rai92}, the author study electric perturbations of a Schr\"odinger operator with periodic potential, whose band functions have non-degenerate minima, providing localization in phase space. The general method described in this work has then been applied with success to numerous models where the band functions share the same properties, such as: a perturbation of a constantly twisted waveguide (\cite{BriKovRaiSoc09}), a half-plane with Neumann boundary condition and constant magnetic field (model \ref{A}) in \cite{BruMirRai13,DomHisSoc13}, and extended to the study of the spectral shift function for a constant magnetic field in a strip (\cite{BriRaiSoc08}). 

Let us describe, as a model example \cite[Theorem 2.3]{BruMirRai13}, analyzing the discrete spectrum for a perturbation of model \ref{A} with a Neumann boundary condition. Recall that the first band function has a unique non-degenerate minimum $\Theta_{0}$, reached for $k=k_{0}$. Assume that $V\in L^{\infty}(\Omega)$ is non-negative potential, and that $V\to0$ as $|(x,y)|\to+\infty$. Then the number of negative eigenvalue of $H_{0}-V$ is, at first order, the same as the number of negative eigenvalues of the effective Hamiltonian 
$$\gh_{eff}:=-\mu\partial_{y}^2-v,$$
acting on $L^{2}(\R)$, where the effective mass is $\mu:=\frac{1}{2}\lambda_{1}''(k_{0})$, and the effective potential is $v(y):=\int_{x\in \R^{+}}V(x,y)u_{1}(x,k_{0})^2\rd x$. 

One of the main ideas behind this result, from \cite{Rai92}, is to exploit the Taylor expansion of the band function near its minimum: $\lambda_{1}(k)= \mu (k-k_{0})^2+o((k-k_{0})^2)$. Therefore, locally, after a translation in $k$ and an inverse Fourier transform, the symbol of $H_{0}$ is in some sense the same as $-\mu \partial_{y}^2$. The effective potential is what remains from $V$ after a localization in frequency near $k_{0}$.

This strategy does not work directly for the same operator with a Dirichlet boundary condition, because the infimum of the first band function is not a minimum anymore, but is the limit of $\lambda_{1}$ at infinity. Still, an effective operator governs the asymptotics. Let us define 
\bel{E:fonctiononde}
\Psi_{x,\xi}(k)=\pi^{-1/4}e^{-i\xi k}e^{-\frac{(x-k)^2}{2}},
\ee
and denote by $P_{x,\xi}$ the projector onto $\mathrm{span}(\Psi_{x,\xi})$. Introduce the anti-Wick integral operator $\cV:L^{2}(\R)\to L^{2}(\R)$ as
\bel{E:Antiwick}
\cV:=\int_{\Omega}V(x,\xi)P_{x,\xi} \rd x \rd \xi.
\ee
This operator is defined through the bilinear form 
\begin{align}
(f,g)\mapsto \langle \cV f,g \rangle:=&\int_{k}\int_{x,\xi}V(x,\xi) (P_{x,\xi}f)(k)\overline{g(k)} \rd k=\langle \cV f,g \rangle
\\
=&\int_{x,\xi}V(x,\xi) \langle \Psi_{x,\xi},f\rangle \langle \Psi_{x,\xi},g\rangle \rd x \rd \xi.
\end{align}
Then, roughly, the operator $\lambda_{1}+\cV$ plays the role of the effective Hamiltonian, where $\lambda_{1}$ stands for the operator of multiplication by $\lambda_{1}$. We refer to \cite[Theorem 2.1]{BruMirRai13} for the precise statement.

This result has been extended to other situations, mainly to the behavior of the spectral shift function at threshold in \cite{BruMir18}, and to other model for which the band functions have a different behavior at its limit, such as the Iwatsuka model. It is on this case that we will focus in Section \ref{S:perturbIwa}. To conclude, we present another model, in which an infinite number of band function accumulate to the bottom of spectrum: the magnetic wire \ref{C}. We analyze the finiteness of the spectrum for perturbation of this Hamiltonian.

\section{Perturbation of the Iwatsuka model}
\label{S:perturbIwa}
In this section we are concerned about extensions of the results from \cite{Shi03,Shi05,Mir16}, counting the number of eigenvalues in the gap of the bands of the spectrum  of the Iwatsuka model B\ref{B1}. In \cite{Mir16}, an effective operator of the same kind as in \cite{BruMirRai13} appears:

Let $0\leq V \in L^{\infty}(\R^{2})$ be a potential such that $\lim_{+\infty}V=0$. Assume that $\cE_{j}^{+}<\cE_{j+1}^{-}$ (this is satisfied when $B^{-}<3B^{+}$). 
Set $A\in (\cE_{j}^{+},\cE_{j+1}^{-})$, and for $\lambda$ small enough, denote by $\cN_{j}(\lambda):=\Tr\one_{(\cE_{j}^{+}+\lambda,A)}$ the number of eigenvalues in the interval $(\cE_{j}^{+}+\lambda,A)$.

Let $\psi_{j}$ be the $j$-th Hermite function, and similarly to \eqref{E:fonctiononde}, define
$$\Psi_{j,x,\xi}(k)=e^{-i\xi k}\psi_{j}(x-a^{-1}(k)).$$
These functions may be seen as a good approximation of the eigenfunctions of $\gh(k)$ (see \eqref{E:ghk}), as $k$ gets large. Then, we define $\cV_{j}$ as in \eqref{E:Antiwick}. 

Then, in the same spirit as \cite{BruMirRai13}, it is stated in \cite[Theorem 2.1]{Mir16} that for all $\delta\in (0,1)$, as $ \to 0$:
$$\Tr \one_{(\cE_{j}^{+},+\infty)}(\lambda_{j}+(1-\delta)\cV_{j})+\cO(1) \leq \cN_{j}(\lambda) \leq  \Tr \one_{(\cE_{j}^{+},+\infty)}(\lambda_{j}+(1+\delta)\cV_{j})+\cO(1),$$
where $\lambda_{j}$ should be understood as the operator of multiplication by $\lambda_{j}$. Note that in this asymptotics, the term $\cO(1)$ depends on $\delta$, and that changing the constant $A$ in the definition of $\cN_{j}$ will only make a difference of a constant term, which can be included in the $\cO(1)$ term. 


Roughly, this theorem can be exploited in two different cases: polynomially decaying potential, and compactly supported potential. In the first cases, we also ask the magnetic field to decay faster that the potential. Let us go into details for the polynomially decaying potential. 

We will assume that  for any pair $(\alpha, \beta) \in \mathbb{Z}_+^2$ ($\mathbb{Z}_+:=\{0,1,2,...\}$),  there exists a positive constant $C_{\alpha,\beta}$ satisfying 
\bel{apr14}
|\partial_x^\beta\partial_y^\alpha V(x,y)|\leq C_{\alpha, \beta}\langle x, y \rangle^{-m-\alpha-\beta} \quad \mbox{ for all} \,\,(x,y)\in \R^2,
\ee
where $m>0$. 
For $\lambda >0$ set 
\bel{apr20}
N_0(\lambda,V):=\frac{1}{2\pi}vol\{(x,y)\in \R^2; V(x,y)> \lambda, x>0\},
\ee 
where $vol$ denotes  the Lebesgue measure in $\R^2$.
Note that if $V$ satisfies  \eqref{apr14}, then $N_0(\lambda,V) = \cO(\lambda^{-\frac{2}{m}})$,  as $\lambda\to0$. 

We also suppose  that  for some positive constants $C$ and  $\lambda_0$
\bel{jul10a}
N_0(\lambda,V)\geq C \lambda^{-2/m}, \quad 0<\lambda<\lambda_0,
\ee
and  that  $N_0(\lambda,V)$  satisfies a homogeneity  condition of the form
\bel{jul10}
\lim_{\epsilon \downarrow 0}\limsup_{\lambda\downarrow 0}\,\lambda^{2/m}\left(N_0(\lambda(1-\epsilon),V)-N_0(\lambda(1+\epsilon),V)\right)=0.
\ee

 Conditions \eqref{apr14}-\eqref{jul10a}-\eqref{jul10} are necessary to have some control of the function $V$ at infinity, from above and below. A typical condition on $V$ to  satisfy \eqref{apr14}-\eqref{jul10a}-\eqref{jul10}   is: $$\lim_{(x,y)\to \infty} \langle x, y\rangle^{m}V(x,y)=\omega\left(\frac{(x,y)}{|(x,y)|}\right),$$ where  $\omega: S^1\to [\epsilon, \infty)$ is smooth and  $\epsilon>0$.

Let us describe \cite[Corollary 2.4]{Mir16}:  assuming that $B$ is a smooth function, with $B_{+}-B=\cO(\langle x \rangle^{-M})$ as $x\to+\infty$, and if $N_{0}(\lambda,V)$ satisfies \eqref{apr14}, \eqref{jul10a} and \eqref{jul10}, with $1<m<M$, then 
$$\cN_{j}(\lambda)\underset{\lambda\to0}{\sim}B_{+}N_{0}(\lambda,V).$$
On the contrary, if $V$ is compactly supported, then, as in \cite{BruMirRai13}, the counting function of eigenvalues is not related to a volume in phase space, see \cite[Corollary 2.3]{Mir16} for the precise result. In particular, even a compactly supported potential can give rise to an infinite number of eigenvalues.

In \cite{MirPof18}, we have extended some of the above results to the behavior of the spectral shift function. Assume that $V$ satisfies
\bel{HVdecay}
|V(x,y)|\leq C\langle x, y \rangle^{-m} \quad \mbox{ for all} \,\,(x,y)\in \R^2,
\ee
with $m>2$. 
Then the diamagnetic inequality  implies that for any real $E < \inf \sigma(H)$ the operator
 $V^{1/2} (H_0 - E)^{-1}$ is Hilbert--Schmidt, and hence the resolvent difference $(H-E)^{-1}-(H_0-E)^{-1}$ is
a trace-class operator. The last property implies that 
there exists a unique function $\xi = \xi(\cdot; H, H_0) \in L^1(\R; (1+E^2)^{-1}dE),$ called the
Spectral Shift Function (SSF) for the operator pair $(H,H_0)$, that satisfies the Lifshits-Kre\u{\i}n trace formula:
$$\mbox{Tr}(f(H)-f(H_0))=\int_\R \xi(E;H,H_0) f'(E)dE, $$
for each $f \in C_0^\infty(\R)$, and vanishes identically in $(-\infty, \inf \sigma(H))$ \cite{Yaf}. Outside the essential spectrum, the spectral shift function and the eigenvalue counting function agree up to a bounded term, and therefore  the SSF could be seen as an extension of 
the counting function to  the whole real line (although we have to be aware that the SSF is only defined as an element of $L^1(\R; (1+E^2)^{-1}dE)$.)

The following result gives an a priori control on the spectral shift function:
\begin{theorem}\label{thssf}
Suppose that  $V\geq0$  satisfies \eqref{HVdecay} and  that $B$  satisfies \eqref{E:PowerlikeB}. 
Set $H_\pm:=H_0\pm V$, then:
\begin{enumerate}
\item\label{thssf1a}  On any compact set $\mathcal{C}\subset \R\setminus\cT$, 
$\sup_{E\in{\mathcal C}}\xi(E;H_\pm,H_0)<\infty,$
i.e.  the SSF is bounded  away from the thresholds. Moreover, 
as $\lambda \downarrow 0$:
\bel{SSFcontrolthresh}
 \xi(\cE_j^{+}-\lambda;H_\pm,H_0)=\cO(\lambda^{-1}+\lambda^{-2/M}),
\quad \xi({\cE}_j^{+}+\lambda;H_{+},H_0)=\cO(\lambda^{-1}), \ee
\bel{SSFboundedthresh}
\xi(\cE_{j}^{+}+\lambda;H_{-},H_0)=\cO(1).
\ee
\item\label{thssf1} The SSF $\xi(\cdot ;H_\pm,H_0)$ is continuous on $\R \setminus \big(\sigma_p(H_\pm)\cup \cT\big)$, where $\sigma_p(H_\pm)$ denote the set of eigenvalues of $H_\pm$. 
\end{enumerate}
\end{theorem}

From  Theorem  \ref{thssf} we know that the only possible  points where  $\xi(E;H_\pm,H_0)$  is unbounded, \Bk are the thresholds $\cT$. The following result shows that it remains bounded if $V$ tends to 0 fast enough: 
\begin{theorem}\label{thVcompact}
Make the same hypotheses that in Theorem \ref{thssf}, with $m>M+2$.  Then,  as $\lambda\downarrow 0$,  $$\xi(\cE_{j}^{+}+\lambda;H_{\pm},H_0)=\cO(1) \ \ \mbox{and} \ \ \xi(\cE_{j}^{+}-\lambda;H_{\pm},H_0)=\cO(1).$$
\end{theorem}
 This theorem covers four different cases. One of these was already stated, under weaker assumption on $V$, in \eqref{SSFboundedthresh}, and is conceptually very different.

As for the counting function, we are able to precise the behavior of the spectral shift function under additional hypotheses on $V$:
\begin{theorem}\label{thssf2}
Make the same hypotheses that in Theorem \ref{thssf}, with $M>m$. Assume also that $N_0(\lambda, V)$ satisfies  \eqref{jul10a}, \eqref{jul10}. Then,  the following asymptotic formula  at the thresholds holds true: 
\bel{asym_ssf}
\xi(\cE_{j}^{+}\pm\lambda;H_\pm,H_0)=\pm B_+\,{N_0(\lambda,V)}(1+o(1)), \quad \lambda \downarrow 0.
\ee
\end{theorem}

\begin{remark}
Note that \eqref{asym_ssf} could be written in the form
\bel{iwa_tam}\xi(\cE_{j}^{+}\pm\lambda;H_\pm,H_0)
=\pm\int\displaylimits_{\displaystyle{\left\{(x,y)\in \R^2;\, V(x,y)>\lambda, \,x>0\right\}}}B(x)\,\rd x \rd y \,(1+o(1)), \quad \lambda \downarrow 0,\ee
which is in accordance with the formula obtained for the function that counts number of  discrete  eigenvalues  of the  perturbed Pauli operator near zero (see \cite{IwaTam98}). This formula 
 shows that, if $B(x)$ converges to $B_+$ fast enough, the asymptotic behavior of $\xi$ depends only on the limit $B_+$. On the contrary, if the convergence of $B(x)$ to $B_+$ is slow in comparison with the decaying rate of $V$, the formula \eqref{iwa_tam} is no longer true, as is shown by  Theorem \ref{thVcompact}, because the term in the right hand side of \eqref{iwa_tam} is unbounded  as $\lambda\downarrow 0$, for example when $V(x,y)=\langle x, y \rangle^{-m}$ \Bk and $m>M+2$.  Therefore, Theorem \ref{thVcompact} is somehow a non-semiclassical result. \end{remark}

 \begin{remark}\label{r6} 
Let $B_+>B_->0$, and consider the magnetic field satisfying $B(x)=B_+$ for $x>0$, and  $B(x)=B_-$  for  $x<0$ (model B2 in the last chapter). Then, for this model  it should be possible to obtain a formula  similar to \eqref{asym_ssf}, using the results on the band functions of \cite{HisPerPofRay16} combined  with the ideas used to analyze the SSF that appear  in \cite{BruMir18}. 
Further, this may suggests that  a result like \eqref{asym_ssf} is  true as long as $B(x)$ converges fast enough to $B_+$.  However, such a result would need the analysis of band function in the general case. \Bk
\end{remark}

Let us describe briefly the strategy if the proofs. The first step is to consider the Birman-Schwinger type operator $T(z)= V^{1/2} (H_0-z)^{-1} V^{1/2}$, for $z\in \C\setminus \R$, and to show that $\lim_{\delta\downarrow 0}T(E+i\delta) = T(E+i0^{+})$ exists in the norm sense, for every $E \in \R\setminus \cT$.

 For $K$ a compact self-adjoint operator, as Section \ref{S:current}, $\dP_{I}(K)$ denotes the spectral projector associated with $K$ on the interval $I$ . For $r>0$, we set
 $$
n_{\pm}(r; K) : = {\rm Rank}\,\dP_{(r,\infty)}(\pm K);
$$
We exploit a useful representation from \cite[Theorem 1.2]{Push97}: For  almost all $E \in \R$, the SSF  $\xi(E;H_\pm,H_0)$ satisfies  
    \bel{pssf}
    \xi(E;H_\pm,H_0) = \pm \frac{1}{\pi}\int_\R n_{\mp}(1;{\rm Re} \,T(E+i 0)+t\, {\rm Im} \, T(E +i0))\frac{dt}{1+t^2}.
    \ee
Estimates near a threshold and continuity of the limit operator $T(E+i0)$ brings Theorem \ref{thssf}, see \cite[Section 5]{MirPof18}.

The precise behavior at thresholds needs more work. As an example, we give the main lines for the behavior of $\xi(\cE_{j}^{+}-\lambda;H_-,H_0)$ as $0<\lambda \ll 1$, when $M>m$ (Theorem \ref{thssf2}), the other cases being similar. The rank of the imaginary part of $T(E+i0)$ being finite and constant between two thresholds, the possible unbounded contribution comes from its real part, as $E$ approaches a threshold. We write it as a singular integral: 
\bel{E:ReTpv}
\Re T(\cE_{j}^{+}-\lambda+i0^{+})=\mathrm{p.v.} \ \ \int_{\R}\frac{G_{j}(k)^{*}G_{j}(k)}{\cE_{j}^{+}-\lambda_{j}(k)-\lambda}\rd k,
\ee
where $G_{j}$ is the operator valued function defined by  
$$G_j(k)f:=\frac{1}{\sqrt{2\pi}} \int_{\R}\int_\R 
e^{-iky}V(x,y)^{1/2}u_j(x,k)f(x,y)\,\rd x\rd y, \quad f \in L^2(\R^2).$$
Using our hypothesis on the decay of the potential, we are able to isolate its most singular contribution, corresponding to large frequencies $I_{\lambda}:=\{ k\in \R,\lambda_{j}(k)\geq \cE_{j}^{+}-\lambda+\epsilon_{\lambda}\}$, where $\epsilon_{\lambda}\in (0,\lambda)$ will be chosen later (notice that in this interval, the integral written in \eqref{E:ReTpv} is no longer a principal value): we introduce 
$$T_{j}(\lambda):= \int_{I_{\lambda}}\frac{G_{j}(k)^{*}G_{j}(k)}{\cE_{j}^{+}-\lambda_{j}(k)-\lambda}\rd k.$$
This operator writes as $S^{*}S$, where $S=Q \one_{I_{\lambda}}|\lambda_{j}-\cE_{j}^{+}+\lambda|^{-1/2}$, and $Q:L^{2}(\R)\to L^{2}(\R^{2})$ is the integral operator with kernel 
$$(2\pi)^{-1/2}V(x,y)^{1/2}e^{iky}u_{j}(x,k).$$
Let ${\cS}_0^m$ be the class of symbols defined in \cite[Section 6.3.2]{MirPof18}. Since $V\in {\cS}_0^m$, we can deduced as in \cite{Shi05} that 
 \bel{shirai}Q^* \, Q=Op^W(w_j), \quad \mbox{with} \ \ w_j \in {\cS}_0^m,\ee
where $Op^W$ is the standard Weyl quantization. Moreover the symbol $w_j$ satisfies
\bel{sep25}
w_{j}(x,y)={V}(a^{-1}(x),-y)+R_1(x,y),  \ \ \mbox{with} \  \ R_1\in {\cS}^{m+1}_0. 
\ee
 To conclude, we need to approximate the number of eigenvalues near 0 of such a pseudo-differential operators. Using our hypotheses on $V$ and  \cite[Theorem 1.3]{DauRo86}, we are able to estimate this number by a volume in phase space, leading to Theorem \ref{thssf2}.

\section{The magnetic wire: example of an accumulation of band functions}
\label{S:magneticwire}
In this section, $(r,\theta,z)$ denotes cylindrical coordinates of $\R^{3}$. Here we focus on the model \ref{C} from Section \ref{S:threshold}, where $A(r,\theta,z)=(0,0,a(r))$, in the special case $a(r)=\log r$. Here, we denote by $\cF$ the Fourier transform with respect to $z$ and $\Phi$ the angular Fourier transform around the $z$ axis.   We have the direct integral decomposition:
$$ \Phi\cF H_{0}\cF^{*}\Phi^{*}:=\sum^{\bigoplus}_{m \geq 0}\int_{k\in \R}^{\bigoplus} \gh_{m}(k) \rd k$$
where the operators
\bel{D:gmk}
\gh_{m}(k):=-\frac{1}{r}\partial_{r}r\partial_{r}+\frac{m^2}{r^2}+(\log r-k)^2
\ee
acts in $L^2_r(\R_+):=L^2(\R_+,r \rd r)$. For all $(m,k)\in \Z\times \R$ the operator $\gh_{m}(k)$ has compact resolvent. 
We denote by $\lambda_{m,j}(k) $, with $j \geq 1$, the $j$-th band functions,
 associated with a normalized eigenvector $u_{m,j}(\cdot,k)$. 

It is known (\cite{Yaf03}) that $k \mapsto \lambda_{m,j}(k)$ is decreasing with  \Bk
$$\lim_{k\rightarrow - \infty} \lambda_{m,j}(k)= + \infty; \qquad \lim_{k\rightarrow + \infty} \lambda_{m,j}(k)= 0, $$
therefore the spectrum of $H_{0}$ is the same as without magnetic field: $\sigma(H_{0})=\sigma(-\Delta)=[0,+\infty)$. 

\begin{figure}[ht]
\begin{center}
\includegraphics[keepaspectratio=true,width=10cm]{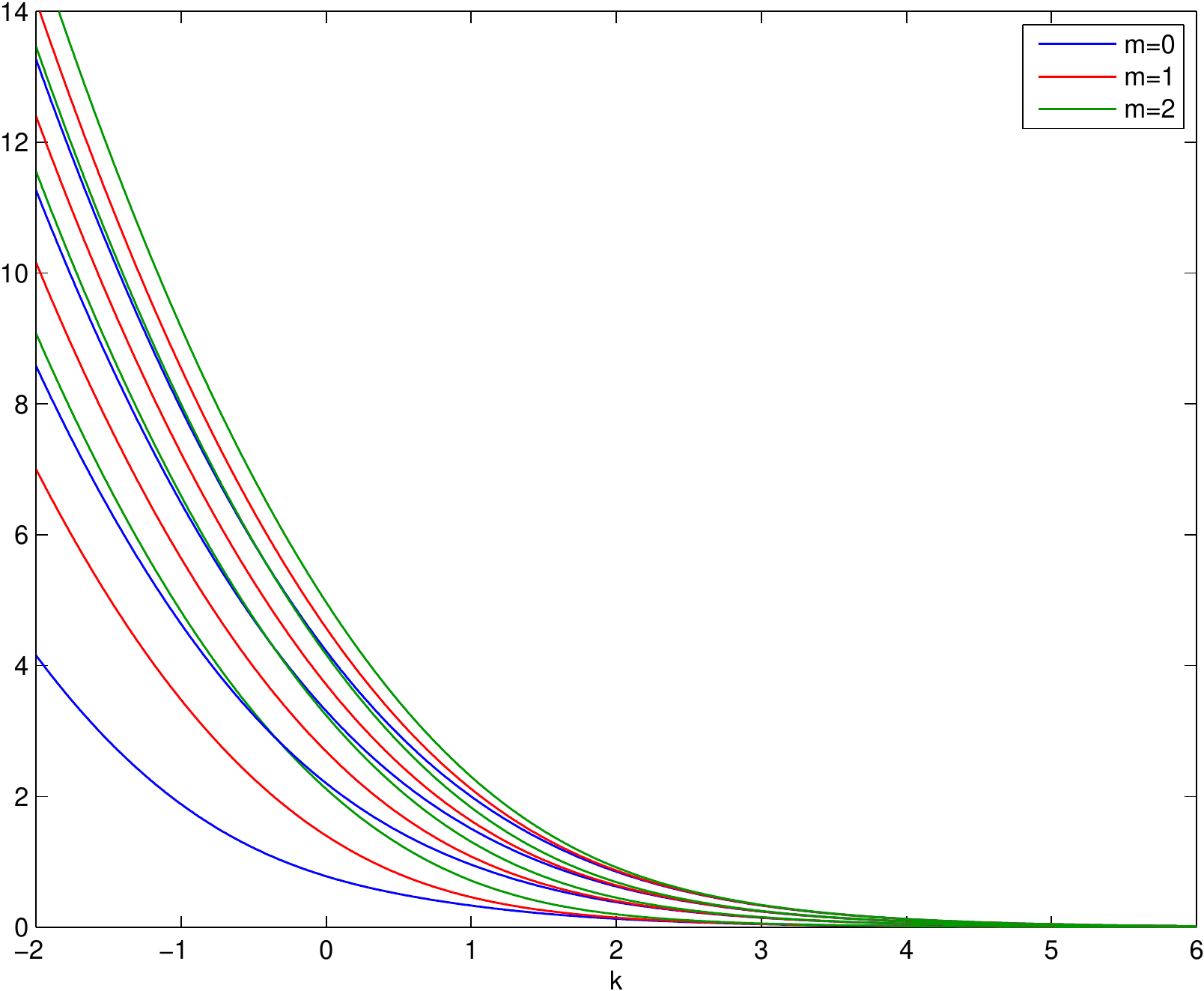}
\caption{The band functions $\lambda_{m,j}(k)$ for $0 \leq m \leq 2$ and $1 \leq j \leq 4$ and $k\in [-2,6]$.}
\label{F1}
\end{center}
\end{figure}

\begin{figure}[h!]
\begin{center}
\includegraphics[keepaspectratio=true,width=10cm]{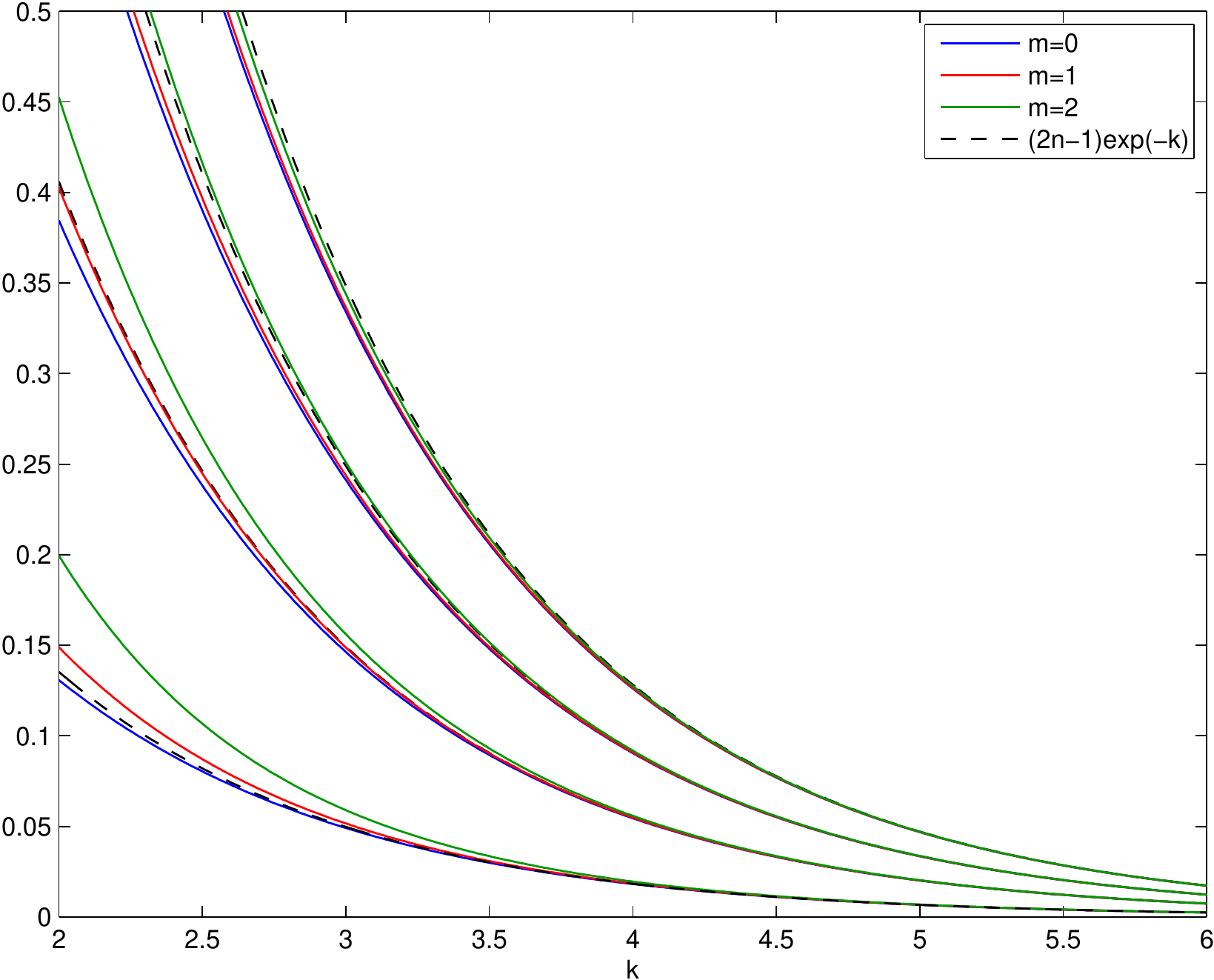}
\caption{Zoom on the lowest energies compared with the first order asymptotics $(2n-1)e^{-k}$. Each cluster corresponds to an energy level $n$.}
\label{F2}
\end{center}
\end{figure}

In \cite[Section 2]{BruPof15} we show that all the band functions tend exponentially to the ground state energy 0 and cluster according to their energy level (see also Figures \ref{F1} and \ref{F2}): 
\begin{theorem}
\label{thm1}
For all $(m,j) \in \Z \times \N^*$, there exist constants $C_{m,j}>0$ and $k_{0} \in \R$ such that for all $k \in (k_{0}, + \infty)$,
\bel{band}
|\lambda_{m,j}(k)-(2j-1)e^{-k}+(m^2-\tfrac{1}{4}-\tfrac{j(j-1)}{2})e^{-2k}| \leq C_{m,j}e^{-5k/2}.
\ee
\end{theorem}
To understand this theorem, as in \cite{Yaf03}, we introduce the parameter 
\bel{E:semicparam}
h:=e^{-k}
\ee
 such that $\log r-k=\log (hr)$.  Some transformations show that $\gh_{m}(k)$ is unitarily equivalent to 
\begin{equation}
\label{D:gh}
\widetilde{\gh_{m}}(h)=-h^2\partial_{\rho}^2+h^2\frac{m^2-\frac{1}{4}}{\rho^2}+\log^2 \rho
\end{equation}
acting on $L^2(\R_{+})$. The positivity of $-\partial_{\rho}^2+\frac{m^2-\frac{1}{4}}{\rho^2}$ involves, with the min-max principle, that the band functions are increasing with respect to $h$, and therefore decreasing with respect to $k$. Moreover, as $h\to0$ (corresponding to $k\to+\infty$), this operator still enters the harmonic approximation, since the leading potential $\log^2\rho$ admits a unique non degenerate minimum, at $\rho=1$. Therefore, the $j$-th eigenvalue of $\widetilde{\gh_{m}}(h)$ will behave as $(2j-1)h$ as $h\to0$. The next term in the asymptotics is of order $h^2$, and involves the parameter $m$ through the Taylor expansion of the whole potential $\log^2\rho+h^2\frac{m^2-\frac{1}{4}}{\rho^2}$ near $\rho=1$, see \eqref{band}. 

We want to know how the ground state energy of $H_{0}$ reacts under electrical perturbation.  For slowly decreasing potentials (with respect to $r$), we have an infinite number of negative eigenvalues of $H_{0}-V$:
\begin{theorem}\label{thm2}
Assume $V$ is a relatively compact perturbation of $H_{0}$ such that
\bel{hypinfinite}
V (x,y,z) \geq r^{-\alpha} \, v_\perp(z), \qquad \alpha>0.
\ee
\Bk
If $\alpha$ and $ v_\perp$ satisfy one of the assumptions {\it (i)}, {\it (ii)} below, then, $H_0-V$ has an infinite number of negative eigenvalues 
which accumulate to $0$.

{\it (i)} $ \alpha<\frac12$ and $v_\perp \in L^1(\R)$ such that
$$\int_\R v_\perp(z)\rd z >0.$$

{\it (ii)} $v_\perp(z) \geq C \langle z \rangle^{-\gamma}$ with $C>0$, $\gamma>0$ and $\alpha+\frac{\gamma}{2} <1$. 
\end{theorem}
This theorem could be expected since the ground state energy of $H_{0}$ corresponds to an infinite number of band functions, therefore it is enough to perturb the free operator by an electric potential which ``reacts'' sufficiently with each band function. The proof relies of the construction of test functions of the form 
$$e^{im\varphi}e^{ikz}u_{m,j}(r,k)f(z), \quad m\in \Z, k\gg 1.$$
Evaluation of these test-functions involves the 1d effective operator $-\partial_{z}^2-V_{m,j}$, where 
$$V_{m,j}(z,k):=\int_{r}|u_{m,j}(r,k)|^2V(r,z)r\rd r.$$
The existence of negative eigenvalues (with a suitable bound) for these operators will provide a negative eigenvalue for $H_{0}-V$, for each $m\in \Z$. This is done through precise estimates of the eigenfunctions $u_{j}(r,k)$ as $k\to+\infty$.

We also have conditions giving finiteness of the negative spectrum. These results are harder to obtain, and sometimes are more surprising.

\begin{theorem}\label{thm3}
Assume $V$ is a relatively compact perturbation of $H_{0}$  such that 
\bel{hypV}
V (x,y,z) \leq \langle (x,y) \rangle^{-\alpha} \, v_\perp(z),
\ee
with $ \alpha>1$, and $v_\perp\in L^p(\R)$ a non negative function with $p\in [1,2]$. 

Then $H_{0}-V$ has, at most, a finite number of negative eigenvalues.
\end{theorem}

 Let us give some comments concerning the above results in comparison with known borderline behavior of perturbations of the Laplacian. Due to the diamagnetic inequality, one might expect for most cases that the density of negative eigenvalues is more important for $-\Delta -V$ than for $H_{0}-V$. Although it is not true in general (see Exemple 2 after Theorem 2.15 of \cite{AvrHerSi78}), the above results illustrate this phenomenon. Theorem \ref{thm2} is a case where the number of negative eigenvalues in the presence of a magnetic field is infinite as without magnetic field. Thanks to Theorem \ref{thm3}, we see that the borderline behavior of the perturbation determining the finiteness of the negative spectrum of $H_{0}-V$ is different than for $-\Delta -V$. \Bk In particular, we obtain:

 


\begin{corollary}\label{cor1}
Let $V$ be a measurable function on $\R^3$ that obeys
$$c \langle (x,y) \rangle^{-\alpha}\langle z \rangle^{-\gamma}\leq V(x,y,z) \leq C \langle (x,y) \rangle^{-\alpha}\langle z \rangle^{-\gamma},$$
with $\alpha+\gamma <2$, $\alpha>1$ and $\gamma >\frac12$.

Then the operator $-\Delta -V$ has infinitely many negative eigenvalues while the negative spectrum of $H_{0}-V$ is finite.
\end{corollary}

Let us sketch the proof of Theorem \ref{thm3}. We introduce the  Birman-Schwinger type operator $T(\lambda)$, decomposed onto the low energies $\{E\leq \nu\}$ and the high energies $\{E > \nu\}$ of $H_{0}$:
\bel{decompT}
T(\lambda)= T_{<\nu}(\lambda) + T_{>\nu}(\lambda),
\ee
with 
$$T_{<\nu}(\lambda) := V^{\frac12} (H_{0}+\lambda)^{-1} {\bf 1}_{[0,\nu]} (H_{0})V^{\frac12}; \qquad 
T_{>\nu}(\lambda) := V^{\frac12} (H_{0}+\lambda)^{-1} {\bf 1}_{]\nu, + \infty[}(H_{0})V^{\frac12}.$$
Here $\nu>0$ is chosen sufficiently small. The Birman-Schwinger principle gives for $\lambda>0$, 
\bel{eq3}
\mathcal{N}(\lambda) = n_+\Big(1, T(\lambda) \Big), 
\ee
The operator  $T_{>\nu}(\lambda)$ is uniformly bounded with respect to $\lambda \geq 0$, therefore we now have to estimate the operator $T_{<\nu}(\lambda)$ corresponding to the low energies. For $V$ axisymmetric, $T_{<\nu}(\lambda)$ is unitarily equivalent to  $\oplus_{m \in \Z} S_m(\lambda)\, S_m(\lambda)^*$ with
 $$S_m(\lambda)\; : \; L^2(\R, l^2(\N^*)) \longrightarrow L^2(\R_+\times \R, \rd r \rd z),$$
 defined, for $(g_j(.))_{j \geq 1} \in  L^2(\R, l^2(\N^*))$ by
\bel{D:Sm}
 S_m(\lambda)(g_j)(r, z):= \frac{V^{\frac12}(r,z)}{\sqrt{2\pi}} \sum_{j \geq 1} \int_{\R} g_j(k) \frac{e^{izk}{\bf 1}_{[0,\nu]}(\lambda_{m,j}(k)) }{(\lambda_{m,j}(k)+\lambda)^\frac12} {\sqrt{r}u_{m,j}(r,k)}\rd k.
 \ee

We estimate the Hilbert-Schmidt norm of this operator: Let $V$ be the axisymmetric potential $V(r,z):=  \langle r \rangle^{-\alpha} \, v_\perp(z)$ with $\alpha>1$ and a non negative function $v_\perp \in L^p(\R)$, $p \in [1,2]$. 
Then  there exist $C>0$ and $\nu_{0}>0$ such that  for all $\nu\in (0,\nu_{0})$ and $\lambda>0$, 
$$\forall m \in \Z, \quad \|S_{m}(\lambda)^{*}S_{m}(\lambda)\|_2 \leq
C \nu^{\alpha-1}.$$
The proof of this estimate relies on a precise estimate of the localization of the eigenfunctions associated with $\lambda_{m,j}(k)$, when $\lambda_{m,j}(k)<\nu$: outside a neighborhood of $e^k$, they decay as $e^{-e^{k}}$, when $k$ gets large, see also \eqref{E:semicparam} for a semi-classical point of view. This is coupled with a uniform lower bound of the band functions such that $\lambda_{m,j}(k)<\nu$, dealing with their behavior as $j\to+\infty$, in order to bound the sum over $j\in \N^{*}$ in \eqref{D:Sm}. The proof of Theorem \ref{thm3} can then be concluded using Weyl's inequality in \eqref{decompT}.


 \bibliographystyle{abbrv}
\bibliography{./bibliopof}

\end{document}